\documentclass{article}
\usepackage{arxiv}

\usepackage{natbib}%

\usepackage{silence}
\usepackage[figuresright]{rotating}
\usepackage{url}
\usepackage{amsthm,amsmath,amsfonts,amssymb}
\usepackage{graphicx} 
\usepackage{subcaption} 

\usepackage{mathtools, nccmath}
\usepackage[ruled,vlined]{algorithm2e}
\newcommand{\lIfElse}[3]{\lIf{#1}{#2 \textbf{else}~#3}}
\usepackage{booktabs,adjustbox}
\usepackage{placeins, paralist, multirow}
\usepackage{float}

\usepackage{pgf,tikz}
\usetikzlibrary{calc}
\usetikzlibrary{arrows,automata,shapes,positioning}

\DeclareMathOperator*{\argmax}{arg\,max}

\SetKw{Break}{break}
\SetKwRepeat{Do}{do}{while}%

\providecommand{\topST}[1]{TOP-ST-MIN#1}
\providecommand{\topP}[1]{TOP-ST-MIN-P#1}
\providecommand{\topPL}[1]{TOP-ST-MIN-PL#1}

\theoremstyle{definition}

\begin{document}
	
	\title{Heuristic approaches for a new variant of the Team Orienteering Problem}
	\author
	{
		Alberto Guastalla\\
		Universit\`a degli Studi di Torino\\
		Dipartimento di Informatica, Torino, Italy\\
		e-mail: \texttt{alberto.guastalla@unito.it} \\
		\and
		Roberto Aringhieri\\
		Universit\`a degli Studi di Torino\\
		Dipartimento di Informatica, Torino, Italy\\
		e-mail: \texttt{roberto.aringhieri@unito.it} \\
		\and
		Pierre Hosteins\\
		Universit\`a degli Studi di Torino\\
		Dipartimento di Informatica, Torino, Italy\\
		e-mail: \texttt{pierre.hosteins@unito.it} \\
	}
	
	\date{\today}
	
	
    \maketitle
    
	
	\begin{abstract}
		In this paper we tackle the Team Orienteering Problem with Service Times, Mandatory Nodes and Incompatibilities, 
		introduced in~\cite{Guastalla2024} and arising from two real-world healthcare applications. We propose two heuristic 
		algorithms in the form of a Variable Descent Neighbourhood algorithm and a matheuristic based on a Cuts Separation 
		approach. For the former, we also provide a multi-thread version exploiting its intrinsic capability to be parallelised. 
		Both algorithms include a specific heuristic routine to provide a starting feasible solution, since finding a feasible 
		solution has been proved to be NP-complete. The results of our heuristic algorithms are compared with an exact cutting 
		plane approach and have complementary strengths and weaknesses. They are also evaluated on existing TOP benchmarks against 
		TOP state-of-the-art algorithms, demonstrating their competitiveness on general grounds.
	\end{abstract}
	
	\keywords{Team orienteering; service times; mandatory nodes; incompatibilities; metaheuritics; matheuristics}

	
	\section{Introduction}
	\label{sec:Introduction}

	The Team Orienteering Problem (TOP) can be described through a complete graph in which each node has a 
	profit and each arc has a travelling time cost. The objective is to find a set of $m$ 
	routes -- starting from the source node and ending at the destination node -- collecting the maximal 
	sum of profits from the visited nodes without exceeding a specific time budget $T_{\max}$ for each of 
	them.

	The TOP has numerous applications (see, e.g., \citet{OP2019}) but the first two in healthcare logistics
	were reported in \citet{ODS}.
	Such applications arise in the logistics of COVID swab collection \citep{DSTCP} and in the ambulance
	routing after a disaster \citep{ARP}, and both share some nice features such as a fixed service
	time, nodes that must be visited, and incompatibilities between nodes meaning that the visit of a node
	is forbidden from certain other nodes.

	In this paper we consider a variant of the TOP called the Team Orienteering Problem with Service Times 
	and Mandatory \& Incompatible Nodes (\topST{}), which has been introduced by \citet{Guastalla2024} as 
	a generalisation of the healthcare applications gathering the three features for the first time in the 
	same TOP model.
	\citeauthor{Guastalla2024} also proved that even finding a feasible solution to the \topST{} is 
	NP-complete, contrary to the TOP.

	In accordance with the literature (see, e.g., \citet{Palomo2015}), we also consider \emph{logical 
		incompatibilities} (usually referred as exclusionary constraints) to provide a more comprehensive 
	analysis.
	Such incompatibilities state that nodes of different types cannot be served by the same route. For 
	instance, in the healthcare setting, patients with an infectious diseases can not share the same 
	vehicles with other not infected patients.
	This led us to consider two versions of the problem: the~\topP{} for the problem with only physical 
	incompatibilities and the~\topPL{} for the problem with both types of incompatibilities.
	
	In \citet{Guastalla2024}, the proposed exact cutting-plane solution algorithm is capable to solve
	small and medium-size instances but is less effective on larger instances.
	In order to find good quality solutions with a reasonable running time for all instances of such 
	a challenging problem, we developed several heuristic algorithms (both for~\topP{} and~\topPL{}) 
	that have been tested on a set of benchmark instances generated in such a way to evaluate the impact 
	of the three features on the problem difficulty.
	We propose two heuristic algorithms and a procedure to quickly find a feasible initial solution.
	The first algorithm is a metaheuristic approach based on the Variable Neighbourhood Descent (VND) methodology. 
	It uses different neighbourhoods in order to improve the overall profit of a solution or to decrease the 
	cost of a route, especially when exploring infeasible solutions.
	The second algorithm consists in a matheuristic approach that dynamically separates the subtour elimination 
	constraints of a linear integer formulation and runs the metaheuristic method as a procedure to fix and improve the current solution. 
	For this algorithm, we also provide a multi-thread version exploiting its intrinsic ability to be parallelised.

	The main contributions of this paper are therefore the following:
	(i) a heuristic solution method to quickly try to find an initial feasible solution, 
	(ii) a metaheuristic approach based on the VND methodology,
	(iii) a matheuristic approach based on the cuts-separation scheme (also presented in a multi-thread version), 
	(iv) an extensive computational analysis including a comparison with the state of the art algorithms for the TOP.

	The paper is organised as follows. 
	Section~\ref{sec:review} reports a survey of the relevant literature to position our contribution among 
	the existing contributions for the TOP and its extensions.
	Section~\ref{sec:algo} provides an accurate description of the proposed solution algorithms.
	Section~\ref{sec:quantitative} presents and discusses the results of the proposed quantitative analysis.
	Section~\ref{sec:end} closes the paper.

	
	\section{Problem statement and literature review}
	\label{sec:review}

	The \topST{} can be modelled on a directed graph $G=(N,A)$ where $N$ is the set of nodes denoted from 
	$1$ (the source or depot) to $n$ (the destination), and $A$ is the set of directed arcs across the 
	nodes in $N$.
	The aim of the \topST{} consists in finding a set of $m$ routes from the source node to the destination 
	one maximising the overall \emph{profit} collected in such a way that the total time (service and travel) of 
	each route does not exceed the maximum allowed time $T_{\max}$, each mandatory node should be visited, 
	and each route should not violate the incompatibilities considered (only physical in the case of
	\topP{} and also the logical ones in the case of \topPL{)}. 
	Hereafter we denote with the \emph{cost of a route} the sum of the travelling and service times of a
	route on a given solution.
	
	
	During the last decades, the literature regarding the orienteering problems has grown very rapidly, 
	underlining their attractiveness for the research community. We focus our attention on the main heuristic 
	methods for the orienteering problems that share at least one of the three features that characterise 
	the \topST{}. 
	In our review, we consider the most recent papers published between the years 2017 and 2023: $11$ are about 
	the Orienteering Problem (OP), while the remaining $19$ are related to the TOP.
	Due to a lack of space we only show the summary Table~\ref{tab:contributions} of the articles considered.
	
	\begin{table}[!ht]  
		\centering  
		\begin{adjustbox}{width=0.89\textwidth,center}  
			\begin{tabular}{@{}llcccccccc@{}}  
				\toprule  
				& \textbf{Reference} & \textbf{Heuristic} & \textbf{Metaheuristic} & \textbf{Matheuristic} & \textbf{(A)} 
				& \textbf{(B)} & \textbf{(C)} & \textbf{(D)} & \textbf{(E)} \\  
				\cmidrule(r){1-2} \cmidrule(lr){3-5} \cmidrule(l){6-10}  
				\textbf{OP}  
				& \citet{Palomo2015} & & \checkmark & & & & \checkmark & & \checkmark \\ 
				& \citet{Wang2017} & \checkmark & \checkmark & & \checkmark & & & \\
				& \citet{ZHENG2017335} & & \checkmark & & \checkmark & & & \\ 
				& \citet{Liao2018} & & \checkmark & & \checkmark & & & \\
				& \citet{Lim2018} & \checkmark & & & \checkmark & & & \\  
				& \citet{LU201854} & & \checkmark & & & & \checkmark & & \checkmark \\
				& \citet{SUN20181058} & \checkmark & & & \checkmark & & & \\     
				& \citet{Exposito2019} & & \checkmark & & \checkmark & & & \\
				& \citet{YU2019488} & & & \checkmark & & \checkmark & & \\  
				& \citet{Dutta2020} & & \checkmark & & \checkmark & & & \\    
				& \citet{Zheng2020} & & \checkmark & & \checkmark & & & \\
				\cmidrule(r){1-2} \cmidrule(lr){3-5} \cmidrule(l){6-10} 
				\textbf{TOP} 
				& \citet{Aghezzaf2015} & & \checkmark & & \checkmark & & & \\ 
				& \citet{Kotiloglu2017} & & \checkmark & & \checkmark & & \checkmark & \\  
				& \citet{LIN2017195} & & \checkmark & & \checkmark & & \checkmark & \\
				& \citet{Sylejmani2017} & & \checkmark & & \checkmark & & & \\
				& \citet{Reyes2018} & & \checkmark & & \checkmark & & & \\
				& \citet{Jin2019} & & \checkmark & & \checkmark & & & \\  
				& \citet{Stavropoulou2019} & & \checkmark & & \checkmark & & \checkmark & \\  
				& \citet{Thompson2019} & & \checkmark & & \checkmark & & & \\  
				& \citet{Gndling2020TimeDependentTT} & & \checkmark & & & \checkmark & & \\ 
				& \citet{Hanafi2020} & & & \checkmark & \checkmark & & & \\  
				& \citet{Wisittipanich2020} & & \checkmark & & \checkmark & & & \\  
				& \citet{Feo2021} & & & \checkmark & \checkmark & & & \\
				& \citet{Heris2021} & & \checkmark & & \checkmark & & & & \\  
				& \citet{Ruiz2021} & & \checkmark & & \checkmark & & \\   
				& \citet{ODS} & & \checkmark & & \checkmark & & \checkmark & \checkmark & \\  
				& \citet{DSTCP} & & \checkmark & & \checkmark & & & & \\  
				& \citet{Correa2024} & & \checkmark & & \checkmark & & \checkmark & & \\  
				& \citet{YU2024121996} & & \checkmark & & \checkmark & & & & \\  
				& \citet{LI2024793} &  &  & \checkmark & \checkmark & \\
				\cmidrule(r){1-2} \cmidrule(lr){3-5} \cmidrule(l){6-10}  
				& \textbf{This paper} & \checkmark & \checkmark & \checkmark & \checkmark & & \checkmark & \checkmark & \checkmark \\  
				\bottomrule  
			\end{tabular}  
		\end{adjustbox}  
		\caption{Summary of recent contributions. The evaluated features are: (A) Fixed Service Times, (B) Variable Service Times, 
			(C) Mandatory Nodes, (D) Physical Incompatibilities, and (E) Logical Incompatibilities.}  
		\label{tab:contributions}  
	\end{table}
	
	The most frequent feature found consists in the presence of service times (fixed and variable) at each node, followed 
	by the mandatory nodes and incompatibilities. Most of the considered works are about the classic Tourist Trip Design 
	Problem (TTDP) in which one or more tourists have to visit a subset of points of interest (POIs). 
	Each POI is characterised by a profit and a service time, and there are several mandatory locations to visit. 
	
	Regarding the solution methods, we mainly found metaheuristic approaches and a few matheuristics and tailored heuristic 
	methods. Table~\ref{tab:contributions} clearly shows that metaheuristics make up for the majority of the solution algorithms
	while matheuristics are less commonly used to solve the TOP.
	
	To the best of our knowledge, the most efficient and effective heuristic solution methods for the TOP are proposed 
	by \citet{HAMMAMI2020105034}. The authors did an extensive literature review collecting the statistics of all the most 
	relevant heuristic algorithms for the TOP and proposed a matheuristic approach based on the hybrid adaptive 
	large neighbourhood search. 
	Among earlier works, \citet{Dang2011} and \citet{DANG2013332} presented two metaheuristic methods based on the particle 
	swarm optimisation scheme, \citet{KIM20133065} reported a metaheuristic approach based on the large neighbourhood 
	search and \citet{VIDAL2014658} proposed a unified solution framework for multi-attribute vehicle routing problems. 
	Finally, \citet{KE2016155} proposed an algorithm called ``Pareto mimic algorithm'' that maintains a population of 
	incumbent solutions which are updated using Pareto dominance.
	All the mentioned algorithms are efficient and effective, albeit with random aspects. 

	\section{Solution algorithms}
	\label{sec:algo}
	
	In this section, we report the three algorithms developed to solve the two versions of the \topST{}
	that we are dealing with, that is the \topP{} and the \topPL{}.
	For the algorithms discussed in this section, the only difference lies in the feasibility test of a 
	given solution, since for the \topPL{} logical incompatibilities should also be checked.
	In the following, we discuss our algorithms directly for the \topPL{}.
	Finally, we observe that finding a feasible solution is NP-complete both for the \topP{} and 
	the \topPL{} since our proof reported in \citet{Guastalla2024} is based on the interaction among
	mandatory nodes and physical incompatibilities.

	The first algorithm (Section~\ref{ssec:maa}) is a procedure which tries to quickly find a feasible initial 
	solution, if it exists.
	The second algorithm (Section~\ref{ssec:metaheuristic}) is a metaheuristic approach based on the variable 
	neighbourhood descent methodology. 
	It uses different neighbourhoods in order to improve the overall profit of a solution or to decrease 
	the cost of a route when exploring infeasible solutions.
	The third algorithm (Section~\ref{ssec:matheuristic}) consists in a matheuristic approach that dynamically 
	separates the subtour elimination constraints and runs the metaheuristic method as a procedure to 
	fix and improve the current solution. 
	
	To simplify the description of the above algorithms, we introduce the set $\hat{N} = \{2, \ldots, n - 1\}$ 
	to represent the set of \emph{customers} and the set 
	$\hat{A} = \{(i,j) \in A : i \in \hat{N} \cup \{1\}, \: j \in \hat{N} \cup \{n\}, \: i \neq j \}$ to represent 
	the set of \emph{traversable} arcs.
	To each node $k \in \hat{N}$ is associated a non-negative amount of service time $s_k$ and a profit $p_k$.
	Each arc $(i,j) \in \hat{A}$ has a non-negative travel time $t_{ij}$ across the node $i$ and $j$.
	Let $c_r$ be the cost of a route $r$, which consists in the sum of $s_k$ and $t_{ij}$ of the nodes and 
	arcs belonging to the route $r$.
	Moreover, we introduce the following three sets: let $M \subset \hat{N}$, $I \subset \hat{A}$, 
	$C \subset \hat{N} \times \hat{N}$ be the set of mandatory nodes, the set of physical incompatibilities, 
	and the set of logical incompatibilities, respectively.

	\subsection{The mandatory nodes allocation algorithm}
	\label{ssec:maa}
	
	As demonstrated in \citet{Guastalla2024}, the feasibility problem of the \topST{} is NP-complete, therefore 
	there is no guarantee to always find a feasible solution for the problem in a short amount of time.
	Given that the difficulty stems from the presence of mandatory nodes, we focus our attention on the set $M$ 
	of mandatory nodes.

	Let $g_k$ be the number of remaining routes in which the mandatory node $k\in M$ could be inserted while 
	respecting the logical incompatibilities: the higher the score, the easiest it is to find how to insert 
	the node in a potential route.
	The Mandatory Nodes Allocation Algorithm (MNAA) starts by defining the set $U$ of unassigned mandatory 
	nodes ($U = M$ at the beginning) setting also $g_u = m$ for each node $u \in U$.
	At each iteration, the MNAA tries to assign a node in $U$ to some route.
	In order to maximise the possibility of finding a feasible solution that includes all the mandatory nodes, 
	the MNAA computes the subset $U_{\overline{g}}$ of unassigned mandatory nodes associated with the minimal score 
	$\overline{g}=\min\{g_u:u\in U\}$. To break ties among the nodes in $U_{\overline{g}}$, the MNAA also computes 
	the next minimal score $\hat{g}_u$ for each node $u \in U_{\overline{g}}$.
	Thus, the MNAA inserts the node $q = \argmax_{u \in U_{\overline{g}}} \hat{g}_u$ in solution, i.e., the one 
	that maximises the next minimal score.
	The insertion of a mandatory node could require also the insertion of further optional nodes to avoid the
	violation of the physical incompatibilities constraints. Let $h^q_{o}$ be the smallest subset of non-mandatory nodes 
	that allow the insertion of $q$, which is then inserted in the solution in such a way to 
	minimise the total service and travelling time.
	After each insertion, the MNAA updates the scores of the unassigned mandatory nodes. 
	If there is at least one mandatory node that cannot be inserted in any route due to physical or logical 
	incompatibilities ($\overline{g} = 0$), the procedure stops returning an infeasibility flag. 
	%
	
	At the end of the MNAA, the solution could still be infeasible because of the violation of the time budget 
	constraint. In our framework, this kind of solution is accepted since the metaheuristic approach, described 
	in the following section, is able to handle it and to regain feasibility.
	
	\begin{algorithm}[!ht]
		\small
		\KwResult{Solution $S$}
		$S \leftarrow \emptyset$ ; $U \leftarrow M$ ; $g \leftarrow m$ ; $g_{\min} \leftarrow m$ \;
		$S \leftarrow$ \texttt{initialiseSolution} ($S$, $m$) \;
		
		\While{$U$ $\textnormal{is}$ $\textnormal{not}$ $\textnormal{empty}$ $\textnormal{and}$ $g_{\min} > 0$}{
			$U_{g_{\min}} \leftarrow$ \texttt{getMinimalScoreNodes} ($g$, $g_{\min}$) \;
			$q \leftarrow$ \texttt{selectMandatoryNode} ($U_{g_{\min}}$) \;
			$h^q_{o}$, $r$, $p$ $\leftarrow$ \texttt{getInsertionList} ($q$, $O$) \;
			$S \leftarrow$ \texttt{updateSolution} ($S$, $h^q_{o}$, $r$, $p$) \;
			$g \leftarrow$ \texttt{updateScores} ($g$) \;
			$g_{\min} \leftarrow$ $\min_{u \in U} g$ \;
			$U \leftarrow U \setminus \{q\}$ \;
		}
		\If{$U$ $\textnormal{is}$ $\textnormal{not}$ $\textnormal{empty}$}{
			\Return{$\emptyset$;}
		}
		\Return{$S$;}
		\caption{\texttt{Mandatory Nodes Allocation Algorithm}} 
		\label{algo:pseudocode-H}
	\end{algorithm}
	
	We report the MNAA pseudocode in Algorithm~\ref{algo:pseudocode-H} that makes use of the following procedures:
	\begin{inparaenum}[(i)]
		\item \texttt{initialiseSolution}: initialises the routes with $m$ empty routes (i.e., with only the source 
		and destination nodes),
		\item \texttt{getMinimalScoreNodes}: calculates the subset of unassigned mandatory nodes associated with the 
		minimal score $\overline{g}$,
		\item \texttt{selectMandatoryNode}: selects the mandatory node to insert into a route that maximises the next 
		minimal score $\overline{g}$ over the subset $U_{\overline{g}}$,
		\item \texttt{getInsertionList}: returns the list of nodes $h^q_{o}$ to insert into the route $r$ in the position $p$,
		\item \texttt{updateSolution}: updates the current solution, and
		\item \texttt{updateScores}: update the scores of each unassigned nodes.
	\end{inparaenum}
	
	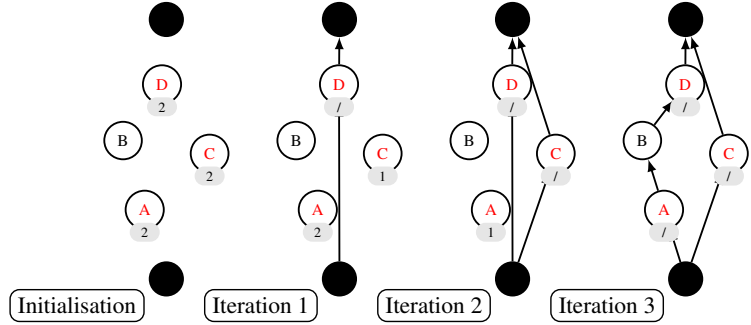
\begin{figure}[!ht]
		\begin{centering}
			\resizebox{0.6\textwidth}{!}{\begin{tikzpicture}[>=latex, scale=0.65, semithick]
    \node [circle, fill, draw, thick, minimum size=0.5cm] (b1) at (0,0) {};
    \node [circle, fill, draw, thick, minimum size=0.5cm] (b2) at (0,6) {};
    \node [circle, fill, draw, thick, minimum size=0.5cm] (b3) at (4,0) {};
    \node [circle, fill, draw, thick, minimum size=0.5cm] (b4) at (4,6) {};
    \node [circle, fill, draw, thick, minimum size=0.5cm] (b5) at (8,0) {};
    \node [circle, fill, draw, thick, minimum size=0.5cm] (b6) at (8,6) {};
    \node [circle, fill, draw, thick, minimum size=0.5cm] (b7) at (12,0) {};
    \node [circle, fill, draw, thick, minimum size=0.5cm] (b8) at (12,6) {};
    
    \node [circle, draw, text=red, thick, minimum size=0.5cm] (A1) at (-0.5,1.6) {\scriptsize A};
    \node [circle, draw, thick, minimum size=0.5cm] (B1) at (-1,3.2) {\scriptsize B};
    \node [circle, draw, text=red, thick, minimum size=0.5cm] (C1) at (1,2.9) {\scriptsize C};
    \node [circle, draw, text=red, thick, minimum size=0.5cm] (D1) at (-0.1,4.5) {\scriptsize D};
    
    \node [circle, draw, text=red, thick, minimum size=0.5cm] (A2) at (3.5,1.6) {\scriptsize A};
    \node [circle, draw, thick, minimum size=0.5cm] (B2) at (3,3.2) {\scriptsize B};
    \node [circle, draw, text=red, thick, minimum size=0.5cm] (C2) at (5,2.9) {\scriptsize C};
    \node [circle, draw, text=red, thick, minimum size=0.5cm] (D2) at (3.99,4.5) {\scriptsize D};
    
    \node [circle, draw, text=red, thick, minimum size=0.5cm] (A3) at (7.5,1.6) {\scriptsize A};
    \node [circle, draw, thick, minimum size=0.5cm] (B3) at (7,3.2) {\scriptsize B};
    \node [circle, draw, text=red, thick, minimum size=0.5cm] (C3) at (9,2.9) {\scriptsize C};
    \node [circle, draw, text=red, thick, minimum size=0.5cm] (D3) at (7.99,4.5) {\scriptsize D};
    
    \node [circle, draw, text=red, thick, minimum size=0.5cm] (A4) at (11.5,1.6) {\scriptsize A};
    \node [circle, draw, thick, minimum size=0.5cm] (B4) at (11,3.2) {\scriptsize B};
    \node [circle, draw, text=red, thick, minimum size=0.5cm] (C4) at (13,2.9) {\scriptsize C};
    \node [circle, draw, text=red, thick, minimum size=0.5cm] (D4) at (11.99,4.5) {\scriptsize D};
    
    \draw [->, thick](b3) to (D2);
    \draw [->, thick](D2) to (b4);

    \draw [->, thick](b5) to (D3);
    \draw [->, thick](D3) to (b6);
    \draw [->, thick](b5) to (C3);
    \draw [->, thick](C3) to (b6);

    \draw [->, thick](b7) to (A4);
    \draw [->, thick](A4) to (B4);
    \draw [->, thick](B4) to (D4);
    \draw [->, thick](D4) to (b8);
    \draw [->, thick](b7) to (C4);
    \draw [->, thick](C4) to (b8);
    
    \node [label, fill=gray!20, rounded corners, minimum width=0.75cm, scale=0.6] at (-0.5, 1.08) {2};
    \node [label, fill=gray!20, rounded corners, minimum width=0.75cm, scale=0.6] at (1, 2.38) {2};
    \node [label, fill=gray!20, rounded corners, minimum width=0.75cm, scale=0.6] at (-0.1, 3.98) {2};
    
    \node [label, fill=gray!20, rounded corners, minimum width=0.75cm, scale=0.6] at (3.5,1.08) {2};
    \node [label, fill=gray!20, rounded corners, minimum width=0.75cm, scale=0.6] at (5, 2.38) {1};
    \node [label, fill=gray!20, rounded corners, minimum width=0.75cm, scale=0.6] at (3.99, 3.98) {/};
    
    \node [label, fill=gray!20, rounded corners, minimum width=0.75cm, scale=0.6] at (7.5,1.08) {1};
    \node [label, fill=gray!20, rounded corners, minimum width=0.75cm, scale=0.6] at (9, 2.38) {/};
    \node [label, fill=gray!20, rounded corners, minimum width=0.75cm, scale=0.6] at (7.99, 3.98) {/};
    
    \node [label, fill=gray!20, rounded corners, minimum width=0.75cm, scale=0.6] at (11.5,1.08) {/};
    \node [label, fill=gray!20, rounded corners, minimum width=0.75cm, scale=0.6] at (13, 2.38) {/};
    \node [label, fill=gray!20, rounded corners, minimum width=0.75cm, scale=0.6] at (11.99, 3.98) {/};

    \node[draw, fill=white, anchor=south east, at={(-0.5, -1)}, rounded corners] {
        Initialisation
    };

    \node[draw, fill=white, anchor=south east, at={(3.5, -1)}, rounded corners] {
        Iteration 1
    };

    \node[draw, fill=white, anchor=south east, at={(7.5, -1)}, rounded corners] {
        Iteration 2
    };

    \node[draw, fill=white, anchor=south east, at={(11.5, -1)}, rounded corners] {
        Iteration 3
    };
\end{tikzpicture}}
			\caption{Algorithm MNAA: small example for problem with $2$ routes.}
			\label{fig:MNAA-plot}
		\end{centering}
	\end{figure}
	
	Figure~\ref{fig:MNAA-plot} depicts how the algorithm works showing the result of the initialisation and three 
	iterations. The figure includes three mandatory nodes (highlighted in red): A, C, and D along with a single 
	optional node B. Note that the value of $\overline{g}$ of a mandatory node is reported in the small grey circle.
	The initialisation assigns the $\overline{g}$ scores to each mandatory node.
	In the first iteration, all nodes have the same score and therefore the MNAA considers the one associated with 
	the minimal insertion cost (in the example, node D). Then, the score of node C is decreased by $1$ since C is 
	logically incompatible with D. 
	In the second iteration, node C is visited by another route because it has the lowest score. Then, the score of 
	node A is then reduced by 1 because it is logically incompatible with C.
	Finally, in the last iteration, node A is visited by the first route through node B due to a physical incompatibility 
	between nodes A and D. This stops its computation since all mandatory nodes have been visited without violating any 
	incompatibility.

	
	\subsection{A variable neighbourhood descent}
	\label{ssec:metaheuristic}
	
	The VND is a deterministic, systematic and efficient local search-based metaheuristic that improves upon the basic local 
	search strategy by exploiting a sequence of different neighbourhood structures \citep{Duarte2018}. The aim is to find a 
	local optimum with respect to one neighbourhood structure before moving on to explore another, potentially leading to 
	better local optima. 
	
	We developed a VND based on two families of neighbourhoods (profit-based and cost-based) that uses a couple of variable-size 
	Tabu-Lists (TLs) to fix some nodes inside or outside the solution for a given amount of iterations.
	In the remainder of this section, we introduce all the ingredients (families of neighbourhoods, search in the feasible 
	and infeasible regions) that are then used to provide the general VND framework.
	
	\subsubsection{Neighbourhoods}
	\label{ssec:neighborhoods}
	
	We consider two families of neighbourhoods with different goals. The profit-based neighbourhoods try to increase the profits of 
	the routes while the cost-based neighbourhoods try to reduce the cost, that is the total travel time of the routes.

	\subsubsection*{The profit-based neighbourhoods}
	
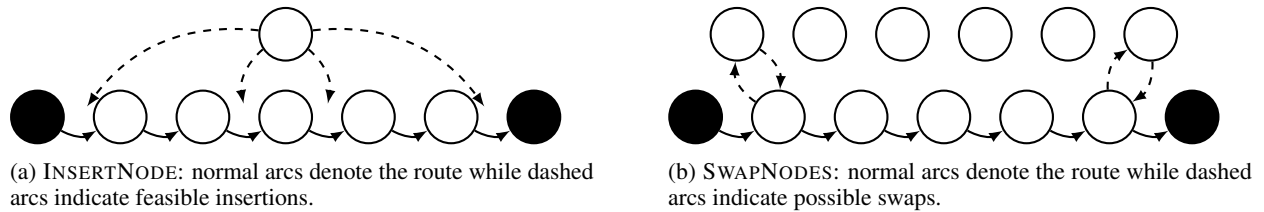
\begin{figure}[!ht]
		\begin{subfigure}[t]{0.47\textwidth}
			\begin{tikzpicture}[>=latex, semithick, scale=0.55]

\node [circle, fill, draw, thick, minimum size=0.7cm] (1) at (0,0) {};
\node [circle, draw, thick, minimum size=0.7cm] (2) at (2,0) {};
\node [circle, draw, thick, minimum size=0.7cm] (3) at (4,0) {};
\node [circle, draw, thick, minimum size=0.7cm] (4) at (6,0) {};
\node [circle, draw, thick, minimum size=0.7cm] (5) at (8,0) {};
\node [circle, draw, thick, minimum size=0.7cm] (6) at (10,0) {};
\node [circle, fill, draw, thick, minimum size=0.7cm] (7) at (12,0) {};
\node [circle, draw, thick, minimum size=0.7cm] (8) at (6,2) {};

\node [circle, draw=none, semithick, minimum size=0.2cm] (f1) at (1,0) {};
\node [circle, draw=none, semithick, minimum size=0.2cm] (f2) at (3,0) {};
\node [circle, draw=none, semithick, minimum size=0.2cm] (f3) at (5,0) {};
\node [circle, draw=none, semithick, minimum size=0.2cm] (f4) at (7,0) {};
\node [circle, draw=none, semithick, minimum size=0.2cm] (f5) at (9,0) {};
\node [circle, draw=none, semithick, minimum size=0.2cm] (f6) at (11,0) {};

\draw [->, thick](1) to [bend right] (2);
\draw [->, thick](2) to [bend right] (3);
\draw [->, thick](3) to [bend right] (4);
\draw [->, thick](4) to [bend right] (5);
\draw [->, thick](5) to [bend right] (6);
\draw [->, thick](6) to [bend right] (7);
\draw [->, thick, dashed](8) to [bend right] (f1);
\draw [->, thick, dashed](8) to [bend right] (f3);
\draw [->, thick, dashed](8) to [bend left] (f4);
\draw [->, thick, dashed](8) to [bend left] (f6);

\end{tikzpicture}
			\caption{\textsc{InsertNode}: normal arcs denote the route while dashed arcs indicate feasible insertions.}
			\label{fig:InsertingNode-plot}
		\end{subfigure}
		\hfill
		\begin{subfigure}[t]{0.47\textwidth}
			\begin{tikzpicture}[>=latex, semithick, scale=0.55]

\node [circle, fill, draw, thick, minimum size=0.7cm] (1) at (0,0) {};
\node [circle, draw, thick, minimum size=0.7cm] (2) at (2,0) {};
\node [circle, draw, thick, minimum size=0.7cm] (3) at (4,0) {};
\node [circle, draw, thick, minimum size=0.7cm] (4) at (6,0) {};
\node [circle, draw, thick, minimum size=0.7cm] (5) at (8,0) {};
\node [circle, draw, thick, minimum size=0.7cm] (6) at (10,0) {};
\node [circle, fill, draw, thick, minimum size=0.7cm] (7) at (12,0) {};

\node [circle, draw, thick, minimum size=0.7cm] (f1) at (1,2) {};
\node [circle, draw, thick, minimum size=0.7cm] (f2) at (3,2) {};
\node [circle, draw, thick, minimum size=0.7cm] (f3) at (5,2) {};
\node [circle, draw, thick, minimum size=0.7cm] (f4) at (7,2) {};
\node [circle, draw, thick, minimum size=0.7cm] (f5) at (9,2) {};
\node [circle, draw, thick, minimum size=0.7cm] (f6) at (11,2) {};

\draw [->, thick](1) to [bend right] (2);
\draw [->, thick](2) to [bend right] (3);
\draw [->, thick](3) to [bend right] (4);
\draw [->, thick](4) to [bend right] (5);
\draw [->, thick](5) to [bend right] (6);
\draw [->, thick](6) to [bend right] (7);
\draw [->, thick, dashed](2) to [bend left] (f1);
\draw [->, thick, dashed](f1) to [bend left] (2);
\draw [->, thick, dashed](6) to [bend left] (f6);
\draw [->, thick, dashed](f6) to [bend left] (6);

\end{tikzpicture}
			\caption{\textsc{SwapNodes}: normal arcs denote the route while dashed arcs indicate possible swaps.}
			\label{fig:SwappingNodeOutside-plot}
		\end{subfigure}
		\caption{Examples of profit-based neighbourhoods.}
		\label{fig:profit-based-neigh}
	\end{figure}

	\begin{compactdesc}
		\item[\textsc{InsertNode}.]
		The aim is to increase the profit of a route, and therefore of the whole solution, by 
		inserting a non assigned node. 
		Such a node is identified by exploring the entire 
		set of insertions then the size of this neighbourhood is $\mathcal{O}\left(n^2\right)$. The best node is 
		then inserted in such a way as to minimise its insertion cost. 
		Figure~\ref{fig:InsertingNode-plot} depicts an example.

		\item[\textsc{SwapNodes}.]
		The goal is to increase the profit by swapping a pair of nodes (one in solution and 
		one not) 
		in a route. This neighbourhood therefore identifies the best swap between
		each pair of nodes and, to break the tie, the swap that yields the smallest cost increase (or even 
		reduces it) is selected. The size of this neighbourhood is $\mathcal{O}\left(n^2\right)$. 
		Figure~\ref{fig:SwappingNodeOutside-plot} depicts an example.
		
	\end{compactdesc}

	\subsubsection*{The cost-based neighbourhoods}
	
\begin{figure}[!ht]
		\begin{subfigure}[t]{0.47\textwidth}
			\begin{tikzpicture}[>=latex, semithick, scale=0.55]

\node [circle, fill, draw, thick, minimum size=0.7cm] (1) at (0,0) {};
\node [circle, draw, thick, minimum size=0.7cm] (2) at (2,0) {};
\node [circle, draw, thick, minimum size=0.7cm] (3) at (4,0) {};
\node [circle, draw, thick, minimum size=0.7cm] (4) at (6,0) {};
\node [circle, draw, thick, minimum size=0.7cm] (5) at (8,0) {};
\node [circle, draw, thick, minimum size=0.7cm] (6) at (10,0) {};
\node [circle, fill, draw, thick, minimum size=0.7cm] (7) at (12,0) {};

\node [circle, draw=none, thick, minimum size=0.1cm] (f1) at (5,0) {};
\node [circle, draw=none, thick, minimum size=0.1cm] (f2) at (7,0) {};

\draw [->, thick](1) to [bend right] (2);
\draw [->, thick](2) to [bend right] (3);
\draw [->, thick](3) to [bend right] (4);
\draw [->, thick](4) to [bend right] (5);
\draw [->, thick](5) to [bend right] (6);
\draw [->, thick](6) to [bend right] (7);
\draw [->, thick, dashed](2) to [bend left=45] (f1);
\draw [->, thick, dashed](6) to [bend right=45] (f2);

\end{tikzpicture}
			\caption{\textsc{MoveNode}: normal arcs denote the route while dashed arcs indicate possible moves.}
			\label{fig:MovingNodeInsideRoute-plot}
		\end{subfigure}
		\hfill
		\begin{subfigure}[t]{0.47\textwidth}
			\begin{tikzpicture}[>=latex, semithick, scale=0.55]

\node [circle, fill, draw, thick, minimum size=0.7cm] (1) at (0,0) {};
\node [circle, draw, thick, minimum size=0.7cm] (2) at (2,0) {};
\node [circle, draw, thick, minimum size=0.7cm] (3) at (4,0) {};
\node [circle, draw, thick, minimum size=0.7cm] (4) at (6,0) {};
\node [circle, draw, thick, minimum size=0.7cm] (5) at (8,0) {};
\node [circle, draw, thick, minimum size=0.7cm] (6) at (10,0) {};
\node [circle, fill, draw, thick, minimum size=0.7cm] (7) at (12,0) {};

\node [circle, fill, draw, thick, minimum size=0.7cm] (f1) at (0,2) {};
\node [circle, draw, thick, minimum size=0.7cm] (f2) at (2,2) {};
\node [circle, draw, thick, minimum size=0.7cm] (f3) at (4,2) {};
\node [circle, draw, thick, minimum size=0.7cm] (f4) at (6,2) {};
\node [circle, draw, thick, minimum size=0.7cm] (f5) at (8,2) {};
\node [circle, draw, thick, minimum size=0.7cm] (f6) at (10,2) {};
\node [circle, fill, draw, thick, minimum size=0.7cm] (f7) at (12,2) {};

\node [circle, draw=none, thick, minimum size=0.1cm] (f8) at (1,0) {};
\node [circle, draw=none, thick, minimum size=0.1cm] (f9) at (3,2) {};
\node [circle, draw=none, thick, minimum size=0.1cm] (f10) at (9,0) {};
\node [circle, draw=none, thick, minimum size=0.1cm] (f11) at (11,2) {};

\draw [->, thick](1) to [bend right] (2);
\draw [->, thick](2) to [bend right] (3);
\draw [->, thick](3) to [bend right] (4);
\draw [->, thick](4) to [bend right] (5);
\draw [->, thick](5) to [bend right] (6);
\draw [->, thick](6) to [bend right] (7);
\draw [->, thick](f1) to [bend left] (f2);
\draw [->, thick](f2) to [bend left] (f3);
\draw [->, thick](f3) to [bend left] (f4);
\draw [->, thick](f4) to [bend left] (f5);
\draw [->, thick](f5) to [bend left] (f6);
\draw [->, thick](f6) to [bend left] (f7);
\draw [->, thick, dashed](f2) to [bend right] (f8);
\draw [->, thick, dashed](2) to [bend right] (f9);
\draw [->, thick, dashed](f6) to [bend right] (f10);
\draw [->, thick, dashed](6) to [bend right] (f11);

\end{tikzpicture}
			\caption{\textsc{MoveNodeRoutes}: normal arcs denote the route while dashed arcs indicate possible moves.}
			\label{fig:MovingNodeBetweenRoutes-plot}
		\end{subfigure}
		\caption{Examples of cost-based neighbourhoods.}
		\label{fig:cost-based-neigh-1}
	\end{figure}
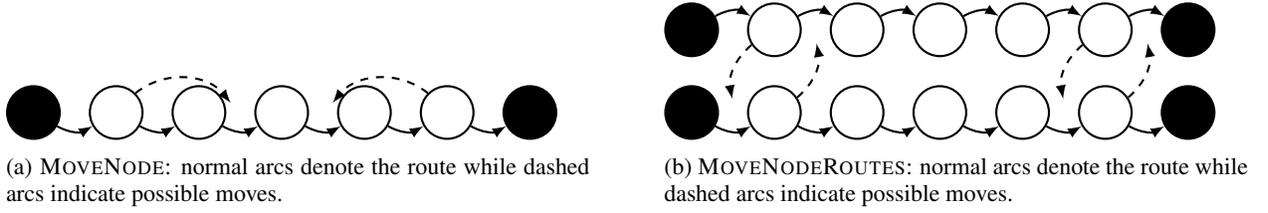
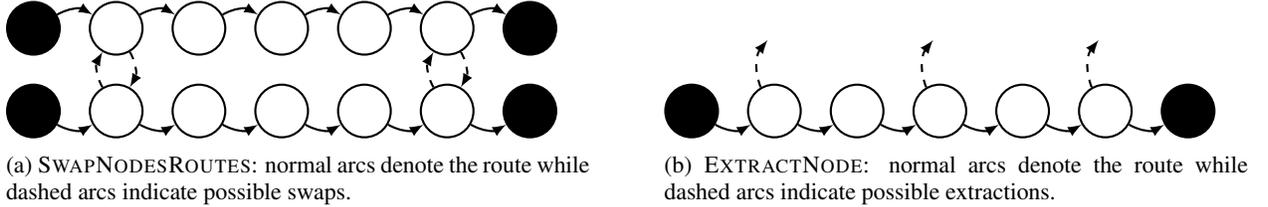

	\begin{compactdesc}
		
		\item[\textsc{MoveNode}.]
		The goal is to reduce the cost of a route by moving a node from one 
		position to another in the same route. To determine the best pair (node, position), that is, the one with 
		the maximal cost reduction, the neighbourhood explores all possible node moves and, therefore, its size is 
		$\mathcal{O}\left(n^2\right)$. 
		Figure~\ref{fig:MovingNodeInsideRoute-plot} depicts an example.

		\item[\textsc{MoveNodeRoutes}.]
		The goal is to reduce the cost by moving a node to another 
		route into a specific position. To determine the best pair (node, position), that is the one with the maximal 
		cost reduction, the neighbourhood explores all possible node moves and, therefore, its size is 
		$\mathcal{O}\left(n^2\right)$. 
		Figure~\ref{fig:MovingNodeBetweenRoutes-plot} depicts an example.
		
	\end{compactdesc}
	
	\begin{compactdesc}
		
		\item[\textsc{SwapNodesRoutes}.]
		The aim of this neighbourhood is to reduce the cost of the solution by swapping a
		pair of nodes belonging to two different routes, maintaining therefore a constant 
		total profit. The best swap is selected by exploring all possible swaps and, therefore, the size of this 
		neighbourhood is $\mathcal{O}\left(n^2\right)$. 
		Figure~\ref{fig:SwappingNodeInside-plot} depicts an example.
		
		\item[\textsc{ExtractNode}.]
		The goal is to reduce the cost of the solution by extracting one node from a route, in such a way 
		to minimise the profit reduction. The size of this neighbourhood is $\mathcal{O}\left(n\right)$. 
		Figure~\ref{fig:ExtractingNodeOutside-plot} depicts an example.
		
	\end{compactdesc}

	\begin{figure}[!ht]
		\begin{subfigure}[t]{0.47\textwidth}
			\begin{tikzpicture}[>=latex, semithick, scale=0.55]
\node [circle, fill, draw, thick, minimum size=0.7cm] (1) at (0,0) {};
\node [circle, draw, thick, minimum size=0.7cm] (2) at (2,0) {};
\node [circle, draw, thick, minimum size=0.7cm] (3) at (4,0) {};
\node [circle, draw, thick, minimum size=0.7cm] (4) at (6,0) {};
\node [circle, draw, thick, minimum size=0.7cm] (5) at (8,0) {};
\node [circle, draw, thick, minimum size=0.7cm] (6) at (10,0) {};
\node [circle, fill, draw, thick, minimum size=0.7cm] (7) at (12,0) {};

\node [circle, fill, draw, thick, minimum size=0.7cm] (f1) at (0,2) {};
\node [circle, draw, thick, minimum size=0.7cm] (f2) at (2,2) {};
\node [circle, draw, thick, minimum size=0.7cm] (f3) at (4,2) {};
\node [circle, draw, thick, minimum size=0.7cm] (f4) at (6,2) {};
\node [circle, draw, thick, minimum size=0.7cm] (f5) at (8,2) {};
\node [circle, draw, thick, minimum size=0.7cm] (f6) at (10,2) {};
\node [circle, fill, draw, thick, minimum size=0.7cm] (f7) at (12,2) {};

\draw [->, thick](1) to [bend right] (2);
\draw [->, thick](2) to [bend right] (3);
\draw [->, thick](3) to [bend right] (4);
\draw [->, thick](4) to [bend right] (5);
\draw [->, thick](5) to [bend right] (6);
\draw [->, thick](6) to [bend right] (7);
\draw [->, thick](f1) to [bend left] (f2);
\draw [->, thick](f2) to [bend left] (f3);
\draw [->, thick](f3) to [bend left] (f4);
\draw [->, thick](f4) to [bend left] (f5);
\draw [->, thick](f5) to [bend left] (f6);
\draw [->, thick](f6) to [bend left] (f7);
\draw [->, thick, dashed](2) to [bend left] (f2);
\draw [->, thick, dashed](f2) to [bend left] (2);
\draw [->, thick, dashed](6) to [bend left] (f6);
\draw [->, thick, dashed](f6) to [bend left] (6);
\end{tikzpicture}
			\caption{\textsc{SwapNodesRoutes}: normal arcs denote the route while dashed arcs indicate possible swaps.}
			\label{fig:SwappingNodeInside-plot}
		\end{subfigure}
		\hfill
		\begin{subfigure}[t]{0.47\textwidth}
			\begin{tikzpicture}[>=latex, semithick, scale=0.55]

\node [circle, fill, draw, thick, minimum size=0.7cm] (1) at (0,0) {};
\node [circle, draw, thick, minimum size=0.7cm] (2) at (2,0) {};
\node [circle, draw, thick, minimum size=0.7cm] (3) at (4,0) {};
\node [circle, draw, thick, minimum size=0.7cm] (4) at (6,0) {};
\node [circle, draw, thick, minimum size=0.7cm] (5) at (8,0) {};
\node [circle, draw, thick, minimum size=0.7cm] (6) at (10,0) {};
\node [circle, fill, draw, thick, minimum size=0.7cm] (7) at (12,0) {};

\node [circle, draw=none, thick, minimum size=0.1cm] (f1) at (2,2) {};
\node [circle, draw=none, thick, minimum size=0.1cm] (f2) at (4,2) {};
\node [circle, draw=none, thick, minimum size=0.1cm] (f3) at (6,2) {};
\node [circle, draw=none, thick, minimum size=0.1cm] (f4) at (8,2) {};
\node [circle, draw=none, thick, minimum size=0.1cm] (f5) at (10,2) {};

\draw [->, thick](1) to [bend right] (2);
\draw [->, thick](2) to [bend right] (3);
\draw [->, thick](3) to [bend right] (4);
\draw [->, thick](4) to [bend right] (5);
\draw [->, thick](5) to [bend right] (6);
\draw [->, thick](6) to [bend right] (7);
\draw [->, thick, dashed](2) to [bend left] (f1);
\draw [->, thick, dashed](4) to [bend left] (f3);
\draw [->, thick, dashed](6) to [bend left] (f5);

\end{tikzpicture}
			\caption{\textsc{ExtractNode}: normal arcs denote the route while dashed arcs indicate possible extractions.}
			\label{fig:ExtractingNodeOutside-plot}
		\end{subfigure}
		\caption{Examples of cost-based neighbourhoods.}
		\label{fig:cost-based-neigh-2}
	\end{figure}
	
	Finally, we also introduced an adaptation of the classical \textsc{2-Opt} neighbourhood \citep{Lin1965} whose
	goal is to reduce the cost of a route by swapping two of its nodes in such a way to maximise 
	the cost reduction. To this end, the full set of swaps are explored and then the size is $\mathcal{O}\left(n^2\right)$.
	
	%
	%
	%
	
	\subsubsection{The exploration of the solution space}
	\label{ssec:exploration-solution-space}
	
	The VND algorithm explores the solution space in two different ways, by alternatively searching for a new solution 
	in the feasible region (SFR phase) or the infeasible region (SIR phase), respectively.
	It makes use of a couple of variable-size Tabu-Lists (TLs) to fix some nodes inside (list $L_1$) or outside 
	(list $L_2$) the solution for a given amount of iterations (list $L_1$ of length $\ell_1$ and list $L_2$ of length 
	$\ell_2$). Their dimension depends on the current solution, and are defined in function of the total number of 
	customers visited by the routes.
	Even if the SFR and SIR phases explore the neighbourhoods in a local search fashion, the two TLs avoid restoring 
	feasibility simply by extracting the node just inserted during the SIR phase. Furthermore, they provide a light 
	intensification and diversification of the search by fixing nodes inside and outside the current solution through 
	the tabu lists.

	The main purpose of the SFR phase is to improve the quality of 
	a solution trying (i) to maximise the profit, and (ii) to minimise the cost while retaining the feasibility 
	of the explored solutions.
	The SFR phase is driven by the following sequence of neighbourhoods: \textsc{InsertNode}, \textsc{SwapNodes}, 
	\textsc{SwapNodesRoutes}, \textsc{MoveNodesRoutes}, \textsc{2-Opt}, and \textsc{MoveNode}. 
	%
	%
	Basically, the idea is to apply all the profit-based neighbourhoods in such a way to improve the profit 
	of the current solution and, then, to apply all the cost-based neighbourhoods to decrease the cost of 
	the obtained routes.
	Each move of the SFR phase is computed by sequentially exploring the above sequence of neighbourhoods 
	until we find the first solution with an better profit or a lower cost (first improvement exploration).
	Each solution $S$ visited during the SFR phase (especially those obtained by applying a cost-based 
	neighbourhood) is improved through a simple intensification phase by greedily adding non-visited 
	nodes while maintaining the cost of the routes as small as possible.
	This intensification is similar to the \textsc{InsertNode} neighbourhood without the use of the
	tabu lists playing therefore a role similar to the \emph{aspiration criteria} of the 
	classical tabu search.
	Finally, the SFR phase continues from the solution $S$ as it allows to explore a larger portion of
	the solution space, as in~\citet{ORP-COR}.

	The main purpose of the SIR phase is to improve the quality of a solution 
	forcing the search to explore an infeasible region of the solution space.
	We define the time-violation of a solution as the sum of the exceeding times of each infeasible route
	$r$, that is the sum of $\max\{c_r - T_{\max},0\}$.
	Therefore, the first step of the SIR phase is to impose the insertion of a node into a route (using the \textsc{InsertNode} 
	neighbourhood) in such a way to make the solution infeasible with respect to the cost, i.e. a solution whose time-violation 
	is strictly greater than zero. 
	The second step consists in exploring the infeasible region in such a way to regain the feasibility by reducing the
	time-violation. To this end, the search is driven by the following sequence of neighbourhoods, that is 
	\textsc{SwapNodesRoutes}, \textsc{MoveNodesRoutes}, \textsc{2-Opt}, \textsc{MoveNode}, \textsc{SwapNodes}, and 
	\textsc{ExtractNode}.
	%
	%
	Each move of the SIR phase is computed by sequentially exploring the above sequence of neighbourhoods 
	until we find the first solution with a lower cost (first improvement exploration).
	It should be noted that the SIR phase cannot guarantee to always restore feasibility, specifically when applying
	moves from the \textsc{ExtractNode} neighbourhood. This could be due to the presence of mandatory nodes and physical 
	incompatibilities in the route whose cost exceeds $T_{\max}$.

	\subsubsection{The algorithm}
	\label{ssec:vnd}
	
	The VND algorithm starts from an initial solution computed by the MNAA, and iterates the SFR and SIR phases.
	At each iteration, the VND will apply a move selected by exploring a certain neighbourhood in accordance with
	the current search phase.
	The VND explores the solution space until a stopping condition is met, that is after reaching a specific number 
	$\zeta_{\text{it}}^{\text{VND}}$ of non-improving iterations (or moves).


	The initial solution could be infeasible since the cost of one or more routes $c_r$ could be greater 
	than $T_{\max}$: in this case, the VND starts its iterations from the SIR phase.
	It may happen that the solution remains infeasible for a number of iterations $\xi$ during the SIR phase. In this case, 
	the VND restarts its computation from a \emph{recovery solution} $S_{\text{R}}$ previously stored, and set to
	the first feasible solution found at the end of the previous SIR phase.
	Finally, if $S_{\text{R}}$ is not available, the VND restarts from the initial solution computed by the MNAA. 
	Note that the exploration of the solution space will be be different from the beginning of the VND due to
	the current status of the TLs whose lengths are updated in accordance with the total number of customers
	visited.

	The VND pseudocode is reported in Algorithm~\ref{algo:pseudocode-MH}, which makes use of the following 
	procedures:
	\begin{inparaenum}[(i)]
		\item \texttt{nextMove}: selects the next move exploring the current (or the next in the list) neighbourhood in the
		current phase ($move = \emptyset$ when all the neighbourhoods in the phase are explored),
		\item \texttt{applyMove}: applies the selected move to the incumbent solution,
		\item \texttt{updateTabuLists}: updates the TLs according to the selected move,
		\item \texttt{resetTabuList}: resets the TLs, and
		\item \texttt{optimise}: inserts the most profitable nodes inside the routes of $S$ minimising their insertion costs.
	\end{inparaenum}
	
	\allowdisplaybreaks{
		\begin{algorithm}[!ht]
			\small
			\KwData{Solution $S$}
			\KwResult{Solution $S^*$}
			$S^* \leftarrow$ $\emptyset$; $S_{\text{R}} \leftarrow \emptyset$ $z^* \leftarrow 0$; $i \leftarrow 0$; $j \leftarrow 0$; 
			$\text{phase} \leftarrow$ (S is feasible ? SFR : SIR);
			$\text{TLs} \leftarrow$ \texttt{resetTabuList}()\;
			\While{not \textnormal{stoppingCondition}($i$)}{
				\If(\tcc*[f]{Search Feasible Region (SFR) phase}){$\textnormal{phase}$ \textnormal{is SFR}} {
					$\text{move}$ $\leftarrow$ \texttt{nextMove} ($S$, $\text{neighbourhoods}$, $\text{TLs}$, $\text{phase}$)\;
					\lIf{$\text{move} = \emptyset$}{
						$\text{phase} \leftarrow$ SIR%
					}
					$j \leftarrow 0$;
				}
				\Else(\tcc*[f]{Search Infeasible Region (SIR) phase}){
					\lIf{$S$ $\textnormal{is}$ $\textnormal{feasible}$}{
						$\text{move}$ $\leftarrow$ \textsc{InsertNode} ($S$, $\text{TLs}$)%
					}
					\lElse{
						$\text{move}$ $\leftarrow$ \texttt{nextMove} ($S$, $\text{neighbourhoods}$, $\text{TLs}$, $\text{phase}$)%
					}
					$j \leftarrow j + 1$;
				}
				
				\If{$\textnormal{move} \neq \emptyset$}{
					$S$ $\leftarrow$ \texttt{applyMove} ($S$, $\text{move}$, $\text{TLs}$)\;
					$\text{TLs}$ $\leftarrow$ \texttt{updateTabuLists} ($\text{move}$, $\text{TLs}$)\;               
				}
				\If{$S$ $\textnormal{is}$ $\textnormal{feasible}$}{
					$\text{phase} \leftarrow$ SFR\;
					$S_{\textnormal{OPT}} \leftarrow $ \texttt{optimise} ($S$)\;
					\lIf{$S_{\text{R}} = \emptyset$}{$S_{\text{R}} \leftarrow S$}
					\lIf{$z^* < z \; (S_{\textnormal{OPT}})$}{
						$S^*$ $\leftarrow S_{\textnormal{OPT}}$; $z^*$ $\leftarrow z \; (S_{\textnormal{OPT}})$%
					}
				}
				\lIf(\tcc*[f]{Restart the search from the recovery solution}){$j = \xi$}{
					$S \leftarrow \textnormal{assignRecoverySolution}(S_{\text{R}})$}         
				$i \leftarrow$ $i + 1$\;
			}        
			\Return{$S^*$;}
			\caption{The \texttt{VND} procedure.} 
			\label{algo:pseudocode-MH}
		\end{algorithm}
	}


	\subsection{A matheuristic algorithm based on a cuts-separation scheme}
	\label{ssec:matheuristic}
	
	The cuts-separation method is one of a variety of optimization methods that iteratively refine a feasible 
	region by finding and adding inequalities to a model.
	A cut can be added to the problem to reduce the search space and this process is repeated until an optimal 
	feasible integer solution is found.
	We propose a Cuts-Separation Heuristic, called CSH.

	We consider a \topST{} formulation based on the extension of the classical formulation for TOP reported in
	\citet{Vansteenwegen2011}, which is characterised by an exponential number  in the number of nodes of 
	subtour elimination constraints (SECs).
	Our algorithm considers a reduced problem (RP) in which all the SECs are removed from the mathematical program.
	Every time a new integer solution for the RP is found, we add all the associated violated SECs to the RP.
	Such a solution is (eventually) repaired and then improved applying the VND.
	To fully exploit the number of feasible solutions generated by the VND, an instance of the Set Packing 
	Problem (SPP) is solved to select the most profitable feasible routes as in \citet{HAMMAMI2020105034}).
	The CSH iterates through the previously described process until a stopping condition is met.

	\subsubsection{The reduced problem}
	\label{ssec:RP}
	
	We formulate the RP on a directed graph defined on $N$, the set of nodes ranging from $1$ (the source or depot) 
	to $n$ (the destination).
	The aim of the RP is to maximise the profit collected by $m$ vehicles without exceeding their time budget $T_{\max}$. 
	Additionally, all the mandatory customers have to be served and no incompatibility between nodes must be violated. 
	The RP can be formulated extending the TOP formulation in~\cite{Vansteenwegen2011} adding the TOP-ST-MIN 
	constraints, and without considering the SECs.
	This formulation makes use of the following decision variables:
	the variable $x_{ijr}$ is equal to $1$ if and only if the arc $(i,j)$ is traversed by the route $r$, $0$ otherwise;
	the variable $y_{kr}$ is equal to $1$ if and only if the node $k$ is visited by the route $r$, $0$ otherwise.
	
	\allowdisplaybreaks{
		\begin{subequations}
			\begin{align}
				\max \quad & {\sum_{r=1}^m \sum_{k=2}^{n-1} p_k \: y_{kr}} \label{eq:RPof}\\
				\textrm{s.t.}
				\quad & \sum_{r=1}^{m} \sum_{j=2}^{n} x_{1jr} = \sum_{r=1}^{m} \sum_{i=1}^{n-1} x_{inr} = m, \label{eq:RPstartend}\\
				\quad & \sum_{i=1}^{n-1} x_{ikr} = \sum_{j=2}^{n} x_{kjr} = y_{kr}, \quad \forall \: k \in \hat{N}, \quad r = 1,\ldots, m, \label{eq:RPconnectivity}\\
				\quad & \sum_{i=1}^{n-1} \sum_{j=2}^{n} (s_i + t_{ij}) \: x_{ijr} \leq T_{\max}, \quad \quad r = 1,\ldots, m, \label{eq:RPbudget}\\
				\quad & \sum_{r=1}^{m} y_{kr} \leq 1, \quad \forall \: k \in \hat{N}, \label{eq:RPonetour}\\
				\quad & \sum_{r=1}^{m} y_{kr} = 1, \quad \forall \: k \in M, \label{eq:RPmandatory}\\
				\quad & \sum_{r=1}^{m} x_{ijr} = 0, \quad \forall \: (i,j) \in I, \label{eq:RPPIncompatibilities}\\
				\quad & |C_k| (1 - y_{kr}) \geq \sum_{i\in {C_k}} y_{ir}, \quad \forall \: k \in \hat{N}, \quad \forall \: r = 1,\ldots, m, \label{eq:RP-logic-incomp}\\
				\quad & x_{ijr} \in \{0, 1\}, \quad \forall \: i = 1, \ldots, n - 1, \quad j = 2, \ldots, n, \quad r = 1,\ldots, m,\label{eq:RPDEF1}\\
				\quad & y_{kr} \in \{0, 1\}, \quad \forall \: k = 2, \ldots, n - 1, \quad r = 1,\ldots, m \label{eq:RPDEF2}
			\end{align}
			\label{mod:RP}
		\end{subequations}
	}
	
	The objective function~\eqref{eq:RPof} maximises the overall profit collected by each route.
	Constraint~\eqref{eq:RPstartend} guarantees (along with the optimisation direction of the 
	objective function) that each route starts from the source node and ends at the destination node.
	Constraints~\eqref{eq:RPconnectivity} guarantee the connectivity of each route while
	constraints~\eqref{eq:RPonetour} guarantee that a node can only be visited by at most one route.
	Constraints~\eqref{eq:RPbudget} guarantee that the total visit time does not exceed $T_{\max}$ for each route.
	Constraints~\eqref{eq:RPmandatory} ensure that all the mandatory customers $M$ will be visited.
	Constraints~\eqref{eq:RPPIncompatibilities} and \eqref{eq:RP-logic-incomp} respectively guarantee that all 
	the physical and logical incompatibilities are satisfied.
	Finally, constraints~\eqref{eq:RPDEF1} and \eqref{eq:RPDEF2} are the variable definition constraints.

	\subsubsection{Repairing and improving the reduced problem solution}
	\label{ssec:fix-solution}
	
	Let $S_{\text{RP}}$ be the solution obtained by solving the RP with some general purpose solver.
	Such a solution can be composed of a set of $m$ routes (from node 1 to node $n$), and some other 
	subtours not connected to the route because of these specific SECs are not included in the RP.

	When such subtours contain some mandatory nodes, the algorithm tries to repair the $m$ routes
	by inserting these mandatory nodes into the routes obtaining a solution $S_M$, which can be
	infeasible due a time-violation. The insertion is performed adopting the same scheme depicted
	in Figure~\ref{fig:MNAA-plot} for the Algorithm~\ref{algo:pseudocode-H}.
	If such an approach is not able to assign all the mandatory nodes, the CSH considers a 
	restricted \topST{} instance composed of the nodes in $S_{\text{RP}}$. Such an instance is
	solved by a general purpose solver exploiting the faster compact formulation reported in
	\citep{Guastalla2024}.
	If the general purpose solver fails to find a feasible solution, the CSH restarts its 
	computation from the recovery solution stored during the VND computation. If this solution
	is not available, the CSH stop returning an infeasibility flag.
	%
	%
	At the end of the repairing routine, the obtained solution $S_{\text{M}}$ will be the initial
	solution for the VND.

	\subsubsection{The set packing problem}
	\label{ssec:SPP}
	
	To further improve the best solution found by the VND, an instances of the SPP is solved to 
	generate an additional feasible solution composed by the most profitable $m$ feasible routes 
	generated in the previous iterations of the VND algorithm.
	In order to guarantee the existence of such a solution, the SPP instance includes all the
	routes generated during the VND whose set of visited customers are different.
	The SPP is formulated as follows:
	
	\allowdisplaybreaks{
		\begin{subequations}
			\begin{align}
				\max \quad & {\sum_r^{\Omega} \lambda_r \: y_r} \label{eq:SPPof}\\
				\textrm{s.t.}
				\quad & \sum_r^{\Omega} a_{kr} y_r \leq 1, \quad \forall \: k \in \hat{N} \setminus M, \label{eq:SPPcustomers}\\
				\quad & \sum_r^{\Omega} a_{kr} y_r = 1, \quad \forall \: k \in M, \label{eq:SPPmandatories}\\
				\quad & \sum_r^{\Omega} y_r \leq m, \quad \forall \: r \in \Omega, \label{eq:SProutes}\\
				\quad & y_r \in \{0, 1\} \quad \forall \: \: r \in \Omega,\label{eq:SPPDef}
			\end{align}
			\label{mod:SPP}
		\end{subequations}
	}
	
	The binary variable $y_r$ is used to select the route $r$ inside the solution.
	The set $\Omega$ represents the whole set of routes computed by several iterations of VND.
	The parameter $a_{kr}$ is 1 if the customer $k$ is visited by the route $r$, 0 otherwise.
	The objective function~\eqref{eq:SPPof} maximises the overall profit of the selected routes.
	Constraint~\eqref{eq:SPPcustomers} guarantees that the non-mandatory customer $k$ is visited 
	by, at most, one route. 
	Constraint~\eqref{eq:SPPmandatories} guarantees that the mandatory customer $k$ is visited by 
	exactly one route. 
	Constraints~\eqref{eq:SProutes} guarantee that a solution is composed by at most $m$ routes.

	\subsubsection{Adding subtour elimination constraints}
	\label{ssec:SEC-addition}
	
	In order to guarantee the progression towards an optimal solution, we add to the RP all the 
	SECs associated with the solution $S_{\text{RP}}$. 
	To this end, we adopted the Depth First Search (DFS) algorithm to compute the set 
	of subtours contained into $S_{\text{RP}}$. In our algorithm, we consider an adaptation of 
	the Dantzig-Fulkerson-Johnson SECs \citep{dfj1954}:
	
	\begin{equation}
		\label{eq:secs}
		\sum_{r=1}^m \sum_{(i,j)\in U \times U} x_{ijr} \leq \sum_{r=1}^m \left(\sum_{i\in U} y_{ir} - y_{kr}\right), 
		\quad \forall \: U \subseteq \{2, \ldots, n-1\}, \; k \in U
	\end{equation}
	
	The validity of~\eqref{eq:secs} follows from the constraints~\eqref{eq:RPconnectivity}
	and~\eqref{eq:RPonetour}.

	\subsubsection{The CSH algorithm}
	\label{algorithm}
	
	A pseudocode of the CSH algorithm is provided in Algorithm~\ref{algo:pseudocode-CP-st}, which
	makes use of the following procedures:
	\begin{inparaenum}[(i)]
		\item \texttt{extractSubtours}: extracts all the subtours from the solution given by the RP,
		\item \texttt{repairSolution}: repairs the RP solution,
		\item \texttt{buildSECs}: builds the SECs associated with the set of subtours, and
		\item \texttt{SPP}: solves an instance of the SPP using the whole set of feasible routes generated.
	\end{inparaenum}
	
	
	\begin{algorithm}[!ht]
		\small
		\KwResult{Solution $S^*$}
		$S^* \leftarrow$ $\emptyset$; $i \leftarrow$ $0$; $\textnormal{SECs}$ $\leftarrow \emptyset$\;
		\For{$i \leftarrow 0$ \KwTo $\zeta_{\text{it}}^{\text{CSH}}$}{
			$S_{\textnormal{RP}} \leftarrow $ \texttt{RP} ($\textnormal{SECs}$)\;
			\If{$S_{\textnormal{RP}} \neq \emptyset$} {
				$\Phi_{\text{RP}} \leftarrow$ \texttt{extractSubtours} ($S_{\textnormal{RP}}$)\;
				\lIfElse{$\Phi_{\text{RP}} \cap M \neq \emptyset$}
				{$S_{\textnormal{M}} \leftarrow$ \texttt{repairSolution} ($S_{\textnormal{RP}}$)}
				{$S_{\textnormal{M}} \leftarrow$ $S_{\textnormal{RP}}$}
				$S_{\textnormal{VND}} \leftarrow$ \texttt{VND} ($S_{\textnormal{M}}$)\;
				$S_{\textnormal{SPP}} \leftarrow$ \texttt{SPP} $(S_{\textnormal{VND}})$\;
				\lIf{$z^* < z \; (S_{\textnormal{SPP}})$}
				{$S^*$ $\leftarrow S_{\textnormal{SPP}}$; $z^*$ $\leftarrow z \; (S_{\textnormal{SPP}})$}
				$\textnormal{SECs}$ $\leftarrow$ $\textnormal{SECs}$ $\cup$ \texttt{buildSECs} ($\Phi_{\textnormal{RP}}$)\;
				$i \leftarrow$ $i + 1$\;
			}
			\lElse{\Return{$\emptyset$}}
		}
		\Return{$S^*$;}
		\caption{The \texttt{CSH-st} single-thread procedure.} 
		\label{algo:pseudocode-CP-st}
	\end{algorithm}

	The CSH starts by initialising the set of SECs to the empty set and the iteration counter.
	The loop is repeated for a maximum number $\zeta_{\text{it}}^{\text{CSH}}$  of iterations.
	Each iteration starts by solving the RP including also the current set of SECs.
	If the RP has no solution, the CSH stops its computation. Otherwise, the CSH
	\begin{inparaenum}[(i)]
		\item extracts the subtours $\Phi_{\text{RP}}$ from the obtained solution $S_{\textnormal{RP}}$,
		\item the solution $S_{\textnormal{RP}}$ is eventually repaired and then improved by applying the VND
		obtaining a new solution $S_{\textnormal{VND}}$,
		\item the $S_{\textnormal{VND}}$ is then refined by solving the corresponding SPP problem, and
		\item the resulting solution $S_{\textnormal{SPP}}$ is finally stored if it improved the current optimal
		solution and the SECs are updated.
	\end{inparaenum}

	From a different perspective, it is worth noting that the CSH generates a sequence of 
	$\zeta_{\text{it}}^{\text{CSH}}$ RP solutions in which eliminating subtours $\Phi_{\text{RP}}$ guides 
	the algorithm towards convergence to an optimal solution.
	The process of improving the RP solutions, i.e., an iteration of Algorithm~\ref{algo:pseudocode-CP-st}, 
	can be parallelised. This leads us to introduce a multi-thread version of the CSH algorithm, which is depicted in
	Algorithm~\ref{algo:pseudocode-CP-mt}:
	\begin{inparaenum}[(i)]
		\item 
		the single-thread part consists in the computation of the $\zeta_{\text{it}}^{\text{CSH}}$ RP solutions
		$S^i_{\textnormal{M}}$, while
		\item 
		the multi-thread part of the algorithm uses a thread pool (whose dimension is set to the number of available threads)
		to solve a $S^i_{\textnormal{M}}$ as soon as a thread in the pool is available.
	\end{inparaenum}
	%

		%
	%
	%
	%
	
	\begin{algorithm}[!ht]
		\small
		\KwResult{Solution $S^*$}
		$S^* \leftarrow \emptyset$; $\textnormal{SECs} \leftarrow \emptyset$\;
		%
		\For(\tcc*[f]{Single-thread execution}){$i \gets 1$ \KwTo $\zeta_{\text{it}}^{\text{CSH}}$}{
			$S_{\textnormal{RP}} \leftarrow $ \texttt{RP} ($\textnormal{SECs}$)
			\If{$S_{\textnormal{RP}} \neq \emptyset$}{
				$\Phi_{\textnormal{RP}} \leftarrow$ \texttt{extractSubtours} ($S_{\textnormal{RP}}$)\;
				\lIfElse{$\Phi_{\text{RP}} \cap M \neq \emptyset$}
				{$S^i_{\textnormal{M}} \leftarrow$ \texttt{repairSolution} ($S_{\textnormal{RP}}$)}
				{$S^i_{\textnormal{M}} \leftarrow$ $S_{\textnormal{RP}}$}
				$\textnormal{SECs}$ $\leftarrow$ $\textnormal{SECs}$ $\cup$ \texttt{buildSECs} ($\Phi_{\textnormal{RP}}$)\;
			}
			\lElse{
				\Return{$\emptyset$}}
			$i \leftarrow$ $i + 1$\;
		}
		$S_{\textnormal{VND}} \leftarrow$ \texttt{VND} ($S^i_{\textnormal{M}}$);
		$S_{\textnormal{SPP}} \leftarrow$ \texttt{SPP} $(S_{\textnormal{VND}})$\tcc*[r]{Multi-thread execution}
		\lIf{$z^* < z \; (S_{\textnormal{SPP}})$}{
			$S^*$ $\leftarrow S_{\textnormal{SPP}}$; 
			$z^*$ $\leftarrow z  \; (S_{\textnormal{SPP}})$}
		\Return{$S^*$;}
		\caption{The \texttt{CSH-mt} multi-thread procedure.} 
		\label{algo:pseudocode-CP-mt}
	\end{algorithm}

	%

	\paragraph{Implementation details}
	
	Two steps of the CSH algorithm could be expensive from a computational point of view,
	that is the solution of the RP and the SPP.
	Here we report the implementation details adopted to limit their running time effort.

	In order to limit the running time to determine the RP solution, we stop the computation of the general purpose
	solver at the root node emphasising the search for feasible solutions fully exploiting the solver heuristics.
	Furthermore, we set a time limit of 360 seconds and a maximal number of 10 feasible solutions found by the
	solver heuristics.

	In order to limit the running time to determine the SPP solution, we introduce an additional stopping condition to the VND 
	that consists in limiting the number of routes generated to 200.000. To consider only different routes, we adopt a hash-table 
	data structure with custom hash function to identify routes with the same set of visited customers.
	Finally, the value of the best VND solution is provided to the solver as a lower bound. 

	
	\section{Quantitative analysis}
	\label{sec:quantitative}
	
	In this section, we report the main findings of an extensive computational analysis to test the VND
	and the CSH algorithm.
	In Section~\ref{ssec:comp-env} we first define the computational environment on which the computational 
	tests are performed.
	Then in Section~\ref{ssec:results-top-st-min} we report the computational results for the \topST{} in 
	both variant \topP{} and \topPL{} comparing their results with those computed by the exact Cutting-Plane 
	Algorithm (CPA) reported in \citet{Guastalla2024}.
	Finally, in Section~\ref{ssec:results-top}  we report also the results on the classical TOP instances 
	obtained by a simple adaptation of our algorithms in order to further validate our approaches.

	\subsection{Setting up the computational environment}
	\label{ssec:comp-env}
	
	The sets of benchmark instances are the ones described in \citet{Guastalla2024}.
	Such instances have been generated in such a way to assess the real impact of the three new features of the
	\topST{}.
	
	The instances were generated starting from the $7$ sets available at the KU Leuven 
	website~\url{https://www.mech.kuleuven.be/en/cib/op} and introduced by \citet{CHAO1996464}.
	Such instances have been grouped by size (\emph{small}, from 21 to 33 nodes, \emph{medium}, from 64 to 66 nodes 
	and \emph{large} from 100 to 102 nodes) and modified in order to include service time, mandatory nodes and 
	physical/logical incompatibilities in such a way that:
	\begin{itemize}
		\item Mandatory nodes: 5\% of the nodes are labelled as mandatory; depending on the instance, those nodes are 
		chosen as \emph{clustered} (close to each other) or \emph{scattered} (far from each other).
		\item Physical incompatibilities: 20\% of the edges are eliminated from the graph, in order to form either 
		\emph{clusters} in the graph (\emph{clusters-based}) or to obtain a \emph{uniform} degree between the nodes 
		(\emph{degree-based}).
		\item Logical incompatibilities: each node is declared logically incompatible with 5\% of the graph nodes, 
		which are either \emph{far} or \emph{near} in terms of travel time.
	\end{itemize}
	These different methods are used for generating the instances are applied in sequence, that is one method for
	the mandatory nodes, then one for the physical incompatibilities, and one for the logical incompatibilities
	only for the \topPL{} instances.
	The acronyms that we will use later to display and analyse the results are summarised in Table~\ref{table:fgs},
	which is divided into two parts: the left side describes the generation methods, while the right side 
	reports the number of instances generated with each method across different instance sizes for both problems.
	Each feature partitions the whole set into two equally sized subsets (e.g., the SM and CM for mandatory nodes 
	divide the set of 156 small instances of \topP{} into two subsets of 78 instances, each containing one half 
	of CPI and one half of DPI instances).
	We refer the interested readers to \citet{Guastalla2024} for additional details on the feature generation
	methods and on the process of instance generation.
	\begin{table}[!ht]
		\begin{adjustbox}{width=1\textwidth,center}
			\begin{tabular}{@{}lllcccccc@{}}
				\toprule
				\multicolumn{3}{c}{\textbf{Generation methods}}                            & \multicolumn{3}{c}{\textbf{TOP-ST-MIN-P}}                                                                     & \multicolumn{3}{c}{\textbf{TOP-ST-MIN-PL}}                                                                    \\
				\textbf{Feature}           & \textbf{Acronym} & \textbf{Description}          & \multicolumn{1}{c}{\textbf{Small}} & \multicolumn{1}{c}{\textbf{Medium}} & \multicolumn{1}{c}{\textbf{Large}} & \multicolumn{1}{c}{\textbf{Small}} & \multicolumn{1}{c}{\textbf{Medium}} & \multicolumn{1}{c}{\textbf{Large}} \\ \cmidrule(r){1-3} \cmidrule(l){4-9}
				Mandatory Nodes            & SM / CM       & Clustered / Scattered         & 78 / 78                                 & 60 / 60                                  & 60 / 60                                 & 156 / 156                                & 120 / 120                                 & 120 / 120                                \\
				Physical Incompatibilities & CPI / DPI     & Clusters-based / Degree-based & 78 / 78                                 & 60 / 60                                  & 60 / 60                                 & 156 / 156                                & 120 / 120                                 & 120 / 120                                \\
				Logical Incompatibilities  & FLI / NLI     & Farthest / Nearest            & --                                 & --                                  & --                                 & 156 / 156                                & 120 / 120                                 & 120 / 120                                \\ \bottomrule
			\end{tabular}
		\end{adjustbox}
		\caption{Summary of feature generation methods for the benchmark instances.}
		\label{table:fgs}
	\end{table}
	
	The VND and the CSH have been implemented in C++ adopting CPLEX 22.1.1 Concert Technology, and
	compiled in Scientific Linux 7.9 (Nitrogen
	The experiments were carried out on a 64-bit Linux machine, with a Six-Core (12 threads) AMD Opteron(tm) Processor 
	8425 HE, 2.10 GHz and 24 GB of RAM. We use CPLEX as a general purpose solver with all built-in cuts in a 
	fully-deterministic mode.
	
	After several preliminary tests, we set the VND parameters as follows:
	\begin{inparaenum}[(i)]
		\item 
		the value of $\zeta_{\text{it}}^{\text{VND}}$ is set to $10,000$ for the small, 
		$100,000$ for the medium, and $500,000$ for the large benchmark set of instances,
		\item 
		the values of $\ell_1$ and $\ell_2$ is set to the $40\%$ of the number of 
		nodes in the current solution (so the length dynamically change),
		\item
		the value of $\xi$ is set to $500$ iterations.
	\end{inparaenum}
	Similarly, we set the CSH parameters as follows:
	\begin{inparaenum}[(i)]
		\item
		the value of $\zeta_{\text{it}}^{\text{CSH}}$ is set equal to $10$ for the small, 
		$30$ for the medium, and $50$ for the large benchmark set of instances,
		\item 
		the value of $\zeta_{\text{it}}^{\text{VND}}$ (for the VNS as inner procedure of CSH) 
		is set equal to $5,000$ for the small, 
		$10,000$ for the medium, and $30,000$ for the large benchmark set of instances,
		\item 
		the CPLEX time limit for the RP and SPP solution is set to $1,200$ seconds.	
	\end{inparaenum}
	
	Finally, our computational results show that \texttt{CSH-mt} provides a speed up, on average, of 
	$7.6$ times the \texttt{CSH-st}. Considering the small, medium and large instances, and
	both \topP{} and \topPL{}, such a speed up ranges between $7.49$ and $7.75$, which represents 
	an honourable value considering that the \texttt{CSH-mt} has an initial not negligible single-thread 
	computation and the maximum possible speed up is $12$ (equal to the maximum number of threads).
	Due to the lack of space, we report only the results for \texttt{CSH-mt} in the remainder 
	of the section, simply denoted by CSH.

	\subsection{Computational results on \topST{} instances}
	\label{ssec:results-top-st-min}
	
	
	Table~\ref{table:VND_CPA_P} reports the comparison between our VND and the CPA on $396$ benchmark 
	instances for the \topP{}, and listed by generation methods as depicted in Table~\ref{table:fgs}.
	The results are summarised in the ALL rows, also by instance size. 
	
	Focusing on the left sub-table,
	\begin{inparaenum}[(i)]
		\item
		the column \# is the number of instances,
		\item
		the column SOL is the number of solutions found,
		\item 
		CPU(s) is the average amount of running time (in seconds), 
		\item 
		LIT(s) is the average last improvement time (in seconds), and
		\item
		LII is the average last improving iteration.
	\end{inparaenum}
	On the right sub-tables, the columns $\Delta_{\text{OBJ}}$ and $\Delta_{\text{CPU}}$ reports 
	the average relative gaps (in percentage) for the values of the objective function and running time, respectively,
	with respect to the CPA.
	Such gaps consider only those instances for which both methods find at least a feasible solution.
	Finally, the column $\Delta_{\text{SOL}}$ reports the absolute gap in terms of number of solutions 
	found.
	
	It is worth noting that $\Delta_{\text{OBJ}}$ and $\Delta_{\text{SOL}}$ are positive when the VND 
	is better than CPA while the value $\Delta_{\text{CPU}}$ is negative when the VND is faster than
	the CPA.
	
	\begin{table}[!ht]
		\centering
		\begin{adjustbox}{width=.75\textwidth,center}
			\begin{tabular}{@{}lrrrrrr@{}}
				\toprule
				\textbf{}       & \textbf{}    & \#  & \textbf{CPU(s)} & \textbf{LIT(s)} & \textbf{LII}    & \textbf{SOL} \\ \cmidrule(r){1-3} \cmidrule(l){4-7}
				\textbf{SMALL}  & \textbf{ALL} & \textbf{156} & \textbf{3.20}    & \textbf{0.01}   & \textbf{66}    & \textbf{156} \\ \cmidrule(lr){2-3} \cmidrule(l){4-7}
				& CM           & 78           & 2.33             & 0.02            & 92             & 78           \\
				& SM           & 78           & 4.07             & 0.01            & 39             & 78           \\
				& CI           & 78           & 2.59             & 0.01            & 57             & 78           \\
				& DI           & 78           & 3.81             & 0.02            & 75             & 78           \\ \cmidrule(r){1-3} \cmidrule(l){4-7}
				\textbf{MEDIUM} & \textbf{ALL} & \textbf{120} & \textbf{100.74}  & \textbf{1.30}   & \textbf{1299}  & \textbf{118} \\ \cmidrule(lr){2-3} \cmidrule(l){4-7}
				& CM           & 60           & 102.00           & 1.51            & 1632           & 60           \\
				& SM           & 60           & 99.47            & 1.08            & 966            & 58           \\
				& CI           & 60           & 111.21           & 1.13            & 954            & 59           \\
				& DI           & 60           & 90.26            & 1.46            & 1644           & 59           \\ \cmidrule(r){1-3} \cmidrule(l){4-7}
				\textbf{LARGE}  & \textbf{ALL} & \textbf{120} & \textbf{1128.27} & \textbf{192.12} & \textbf{97075} & \textbf{94}  \\ \cmidrule(lr){2-3} \cmidrule(l){4-7}
				& CM           & 60           & 1390.93          & 288.62          & 154369         & 60           \\
				& SM           & 60           & 865.62           & 95.63           & 39780          & 34           \\
				& CI           & 60           & 1159.04          & 177.67          & 75371          & 47           \\
				& DI           & 60           & 1097.51          & 206.57          & 118778         & 47           \\ \cmidrule(r){1-3} \cmidrule(l){4-7}
				\textbf{ALL}    &              & \textbf{396} & \textbf{373.69}  & \textbf{58.62}  & \textbf{29836} & \textbf{368}  \\ \bottomrule
			\end{tabular}
			\hspace{3mm}
			\begin{tabular}{@{}lrrrrr@{}}
				\toprule
				\textbf{}       &              & \textbf{$\Delta_{\text{OBJ}}$} & \textbf{$\Delta_{\text{CPU}}$} & \textbf{$\Delta_{\text{SOL}}$} \\ \cmidrule(r){1-2} \cmidrule(l){3-5}
				\textbf{SMALL}  & \textbf{ALL} & \textbf{-5.27} & \textbf{-81.71} & \textbf{0}   \\ \cmidrule(r){2-2} \cmidrule(l){3-5}
				& CM           & -5.82          & -92.13          & 0            \\
				& SM           & -4.64          & -24.75          & 0            \\
				& CI           & -5.36          & -83.50          & 0            \\
				& DI           & -5.17          & -80.25          & 0            \\ \cmidrule(r){1-2} \cmidrule(l){3-5}
				\textbf{MEDIUM} & \textbf{ALL} & \textbf{-3.32} & \textbf{-92.31} & \textbf{1}   \\ \cmidrule(lr){2-2} \cmidrule(l){3-5} 
				& CM           & -2.59          & -93.20          & 0            \\
				& SM           & -4.10          & -91.12          & 1            \\
				& CI           & -2.74          & -90.25          & 0            \\
				& DI           & -3.92          & -93.89          & 1            \\ \cmidrule(r){1-2} \cmidrule(l){3-5}
				\textbf{LARGE}  & \textbf{ALL} & \textbf{-1.79} & \textbf{-81.30} & \textbf{47}  \\ \cmidrule(lr){2-2} \cmidrule(l){3-5} 
				& CM           & -1.42          & -80.43          & 19           \\
				& SM           & -4.25          & -82.56          & 28           \\
				& CI           & -3.10          & -81.07          & 25           \\
				& DI           & -0.61          & -81.55          & 22           \\ \cmidrule(r){1-2} \cmidrule(l){3-5}
				\textbf{ALL}    &              & \textbf{-3.45} & \textbf{-83.26} & \textbf{48} \\ \bottomrule
			\end{tabular}
		\end{adjustbox}
		\caption{The VND results (left) and the comparison with the CPA for the \topP{} (right).}
		\label{table:VND_CPA_P}
	\end{table}
	
	The results reported in Table~\ref{table:VND_CPA_P} show that the VND is always much 
	faster than the CPA computing good quality solutions in terms of objective function.
	A clear advantage of the VND is to find almost always at least one feasible solution, even for those instances that the CPA 
	is not able to certify as infeasible. As a matter of fact, the VND is able to find 48 additional solutions for the 
	\topP{,} 47 of which are solutions for the larger instances.
	
	Table~\ref{table:CSH_CPA_P} reports the comparison between our CSH and the CPA on $396$ benchmark 
	instances for the \topP{}, and listed by generation methods as depicted in Table~\ref{table:fgs}.
	Table~\ref{table:CSH_CPA_P} has the same structure of Table~\ref{table:VND_CPA_P} except for the columns 
	LIT(s) and LII, which are specific for the VND.
	We recall that $\Delta_{\text{OBJ}}$ and $\Delta_{\text{SOL}}$ are positive when the CSH 
	is better than CPA while the value $\Delta_{\text{CPU}}$ is negative when the CSH is faster than
	the CPA.
	
	\begin{table}[!ht]
		\centering
		\begin{adjustbox}{width=.60\textwidth,center}
			\begin{tabular}{@{}lrrrrrr@{}}
				\toprule
				\textbf{}       & \textbf{}    & \# & \textbf{CPU(s)}  & \textbf{SOL} \\ \cmidrule(r){1-3} \cmidrule(l){4-5}
				\textbf{SMALL}  & \textbf{ALL} & \textbf{156} & \textbf{86.62}  & \textbf{156}  \\ \cmidrule(r){2-3} \cmidrule(l){4-5}
				& CM           & 78           & 89.37           & 78            \\
				& SM           & 78           & 75.87           & 78            \\
				& CI           & 78           & 89.07           & 78            \\
				& DI           & 78           & 76.18           & 78            \\ \cmidrule(r){1-3} \cmidrule(l){4-5}
				\textbf{MEDIUM} & \textbf{ALL} & \textbf{120} & \textbf{486.02} & \textbf{118}  \\ \cmidrule(r){2-3} \cmidrule(l){4-5}
				& CM           & 60           & 525.71          & 60            \\
				& SM           & 60           & 446.33          & 58            \\
				& CI           & 60           & 475.88          & 59            \\
				& DI           & 60           & 496.16          & 59            \\ \cmidrule(r){1-3} \cmidrule(l){4-5}
				\textbf{LARGE}  & \textbf{ALL} & \textbf{120} & \textbf{2340.99} & \textbf{94}   \\ \cmidrule(r){2-3} \cmidrule(l){4-5}
				& CM           & 60           & 2556.64          & 60            \\
				& SM           & 60           & 2125.34          & 34            \\
				& CI           & 60           & 2433.48          & 47            \\
				& DI           & 60           & 2248.50          & 47            \\ \cmidrule(r){1-3} \cmidrule(l){4-5}
				\textbf{ALL}    &              & \textbf{396} & \textbf{889.22} & \textbf{368}   \\ \bottomrule
			\end{tabular}
			\hspace{3mm}
			\begin{tabular}{@{}lrrrrr@{}}
				\toprule
				\textbf{}       &              & \textbf{$\Delta_{\text{OBJ}}$} & \textbf{$\Delta_{\text{CPU}}$} & \textbf{$\Delta_{\text{SOL}}$} \\ \cmidrule(r){1-2} \cmidrule(l){3-5}
				\textbf{SMALL}  & \textbf{ALL} & \textbf{-1.06} & \textbf{372.23} & \textbf{0}   \\ \cmidrule(r){2-2} \cmidrule(l){3-5}
				& CM           & -0.97          & 202.12           & 0            \\
				& SM           & -1.16          & 1302.31          & 0            \\
				& CI           & -1.27          & 466.89          & 0            \\
				& DI           & -0.85          & 295.10          & 0            \\ \cmidrule(r){1-2} \cmidrule(l){3-5}
				\textbf{MEDIUM} & \textbf{ALL} & \textbf{-0.06} & \textbf{-62.89}   & \textbf{1}   \\ \cmidrule(r){2-2} \cmidrule(l){3-5}
				& CM           & 0.02           & -64.93             & 0            \\
				& SM           & -0.15          & -60.15            & 1            \\
				& CI           & -0.01          & -58.29            & 0            \\
				& DI           & -0.11          & -66.44            & 1            \\ \cmidrule(r){1-2} \cmidrule(l){3-5}
				\textbf{LARGE}  & \textbf{ALL} & \textbf{1.08}  & \textbf{-61.21}   & \textbf{47}  \\ \cmidrule(r){2-2} \cmidrule(l){3-5}
				& CM           & 1.08           & -64.03             & 19           \\
				& SM           & 1.06           & -57.17            & 28           \\
				& CI           & 0.59           & -60.25            & 25           \\
				& DI           & 1.51           & -62.20            & 22           \\ \cmidrule(r){1-2} \cmidrule(l){3-5}
				\textbf{ALL}    &              & \textbf{-0.05} & \textbf{-60.17}   & \textbf{48}  \\ \bottomrule
			\end{tabular}
		\end{adjustbox}
		\caption{The CSH results (left) and the comparison with the CPA for the \topP{} (right).}
		\label{table:CSH_CPA_P}
	\end{table}
	
	The results reported in Table~\ref{table:CSH_CPA_P} shows that the CSH is capable to compute high quality solutions for all
	the small, medium, and large instances, on average. Regarding the running time, the CSH performs worse on the small instances 
	and better on the medium and large instances. 
	Even in this case, the matheuristic approach is able to find 48 additional solutions for the \topP{,} 47 of which are solutions 
	for the larger instances, similar to the VND.
	
	Tables~\ref{table:VND_CPA_P} and~\ref{table:CSH_CPA_P} reveal some common behaviour regarding the type of instances.
	It is worth noting that the VND and the CSH find more feasible solutions for the SM instances compared to the other 
	ones. Specifically, they identified $28$ additional feasible solutions for the SM instances, compared to 19 for CM 
	instances. On the contrary, the comparison between CI and DI instances revealed no notable differences in the number of 
	feasible solutions found.
	
	
	Tables~\ref{table:VND_CPA_PL} and~\ref{table:CSH_CPA_PL} report the results of our algorithms on the $792$ \topPL{} 
	instances, and they has the same structure of Table~\ref{table:VND_CPA_P} and~\ref{table:CSH_CPA_P}, respectively.
	We recall that $\Delta_{\text{OBJ}}$ and $\Delta_{\text{SOL}}$ are positive when the VND or CSH 
	is better than CPA while the value $\Delta_{\text{CPU}}$ is negative when the VND or CSH is faster than
	the CPA.
	
	\begin{table}[!ht]
		\centering
		\begin{adjustbox}{width=0.70\textwidth,center}
			\begin{tabular}{@{}lrrrrrr@{}}
				\toprule
				\textbf{}       & \textbf{}    & \# & \textbf{CPU(s)}  & \textbf{LIT(s)} & \textbf{LII}      & \textbf{SOL} \\ \cmidrule(r){1-3} \cmidrule(l){4-7}
				\textbf{SMALL}  & \textbf{ALL} & \textbf{312} & \textbf{4.05}    & \textbf{0.02}   & \textbf{65}    & \textbf{312} \\ \cmidrule(r){2-3} \cmidrule(l){4-7}
				& CM           & 156          & 3.23             & 0.02            & 91             & 156          \\
				& SM           & 156          & 4.87             & 0.01            & 40             & 156          \\
				& CI           & 156          & 3.49             & 0.02            & 61             & 156          \\
				& DI           & 156          & 4.60             & 0.02            & 69             & 156          \\
				& FI           & 156          & 4.99             & 0.01            & 53             & 156          \\
				& NI           & 156          & 3.11             & 0.02            & 78             & 156          \\ \cmidrule(r){1-3} \cmidrule(l){4-7}
				\textbf{MEDIUM} & \textbf{ALL} & \textbf{240} & \textbf{163.74}  & \textbf{3.74}   & \textbf{2632}  & \textbf{236} \\ \cmidrule(r){2-3} \cmidrule(l){4-7}
				& CM           & 120          & 173.20           & 5.70            & 4131           & 120          \\
				& SM           & 120          & 154.28           & 1.79            & 1133           & 116          \\
				& CI           & 120          & 169.08           & 2.42            & 1417           & 118          \\
				& DI           & 120          & 158.40           & 5.07            & 
				3847          & 118          \\
				& FI           & 120          & 188.08           & 2.40            & 1177           & 118          \\
				& NI           & 120          & 139.41           & 5.09            & 4087           & 118          \\ \cmidrule(r){1-3} \cmidrule(l){4-7}
				\textbf{LARGE}  & \textbf{ALL} & \textbf{240} & \textbf{1638.88} & \textbf{282.98} & \textbf{59971} & \textbf{142} \\ \cmidrule(r){2-3} \cmidrule(l){4-7}
				& CM           & 120          & 1962.45          & 420.25          & 87648          & 80           \\
				& SM           & 120          & 1315.31          & 145.72          & 32293          & 62           \\
				& CI           & 120          & 1718.79          & 310.02          & 59053          & 71           \\
				& DI           & 120          & 1558.97          & 255.95          & 60888          & 71           \\
				& FI           & 120          & 2575.79          & 502.96          & 91739          & 94           \\
				& NI           & 120          & 701.97           & 63.00           & 28202          & 48           \\ \cmidrule(r){1-3} \cmidrule(l){4-7}
				\textbf{ALL}    &              & \textbf{792} & \textbf{547.84}  & \textbf{86.89}  & \textbf{18996} & \textbf{690} \\ \bottomrule
			\end{tabular}
			\hspace{3mm}
			\begin{tabular}{@{}lrrrrr@{}}
				\toprule
				\textbf{}       &              & \textbf{$\Delta_{\text{OBJ}}$} & \textbf{$\Delta_{\text{CPU}}$} & \textbf{$\Delta_{\text{SOL}}$} \\ \cmidrule(r){1-2} \cmidrule(l){3-5}
				\textbf{SMALL}  & \textbf{ALL} & \textbf{-4.00} & \textbf{-95.81} & \textbf{0}   \\ \cmidrule(r){2-2} \cmidrule(l){3-5}
				& CM           & -4.44          & -98.08          & 0            \\
				& SM           & -3.50          & -80.80          & 0            \\
				& CI           & -4.10          & -96.81          & 0            \\
				& DI           & -3.90          & -94.51          & 0            \\
				& FI           & -5.14          & -90.94          & 0            \\
				& NI           & -2.69          & -97.75          & 0            \\ \cmidrule(r){1-2} \cmidrule(l){3-5}
				\textbf{MEDIUM} & \textbf{ALL} & \textbf{-2.07} & \textbf{-95.95} & \textbf{9}   \\ \cmidrule(r){2-2} \cmidrule(l){3-5}
				& CM           & -1.49          & -96.29          & 0            \\
				& SM           & -2.77          & -95.48          & 9            \\
				& CI           & -2.05          & -96.11          & 5            \\
				& DI           & -2.10          & -95.76          & 4            \\
				& FI           & -2.12          & -94.04          & 5            \\
				& NI           & -2.01          & -97.17          & 4            \\ \cmidrule(r){1-2} \cmidrule(l){3-5}
				\textbf{LARGE}  & \textbf{ALL} & \textbf{4.05}  & \textbf{-72.20} & \textbf{90}  \\ \cmidrule(r){2-2} \cmidrule(l){3-5}
				& CM           & 4.09           & -72.66          & 30           \\
				& SM           & 2.41           & -71.49          & 60           \\
				& CI           & 3.04           & -70.88          & 45           \\
				& DI           & 5.06           & -73.54          & 45           \\
				& FI           & 3.94           & -57.56          & 49           \\
				& NI           & 5.17           & -87.73          & 41           \\ \cmidrule(r){1-2} \cmidrule(l){3-5}
				\textbf{ALL}    &              & \textbf{-1.76} & \textbf{-82.04} & \textbf{99} \\ \bottomrule
			\end{tabular}
		\end{adjustbox}
		\caption{The VND results (left) and the comparison with the CPA for the \topPL{} (right).}
		\label{table:VND_CPA_PL}
	\end{table}
	
	The results reported in Table~\ref{table:VND_CPA_PL} show that the VND is always much 
	faster than the CPA
	, on average.
	Contrary to what is shown on \topP{} instances, the VND computes better quality solutions for the larger
	instances. Furthermore, the VND is capable to find 99 additional solutions with respect to the CPA,
	90 of which are solutions for the larger instances.
	
	\begin{table}[!ht]
		\centering
		\begin{adjustbox}{width=0.60\textwidth,center}
			\begin{tabular}{@{}lrrrrrr@{}}
				\toprule
				\textbf{}       & \textbf{}    & \# & \textbf{CPU(s)}  & \textbf{SOL} \\ \cmidrule(r){1-3} \cmidrule(l){4-5}
				\textbf{SMALL}  & \textbf{ALL} & \textbf{312} & \textbf{82.03}  & \textbf{312} \\ \cmidrule(r){2-3} \cmidrule(l){4-5}
				& CM           & 156          & 101.46           & 156          \\
				& SM           & 156          & 62.60           & 156          \\
				& CI           & 156          & 82.84           & 156          \\
				& DI           & 156          & 81.22           & 156          \\
				& FI           & 156          & 80.26           & 156          \\
				& NI           & 156          & 83.80           & 156          \\ \cmidrule(r){1-3} \cmidrule(l){4-5}
				\textbf{MEDIUM} & \textbf{ALL} & \textbf{240} & \textbf{584.53} & \textbf{236} \\ \cmidrule(r){2-3} \cmidrule(l){4-5}
				& CM           & 120          & 623.79          & 120          \\
				& SM           & 120          & 545.26          & 116          \\
				& CI           & 120          & 594.72          & 118          \\
				& DI           & 120          & 574.33          & 118          \\
				& FI           & 120          & 476.16          & 118          \\
				& NI           & 120          & 692.90          & 118          \\
				\cmidrule(r){1-3} \cmidrule(l){4-5}
				\textbf{LARGE}  & \textbf{ALL} & \textbf{240} & \textbf{2369.71} & \textbf{142} \\ \cmidrule(r){2-3} \cmidrule(l){4-5}
				& CM           & 120          & 2483.53          & 80          \\
				& SM           & 120          & 2255.89          & 62           \\
				& CI           & 120          & 2447.66          & 71           \\
				& DI           & 120          & 2291.77          & 71           \\
				& FI           & 120          & 3500.85         & 94           \\
				& NI           & 120          & 1238.58          & 48           \\
				\cmidrule(r){1-3} \cmidrule(l){4-5}
				\textbf{ALL}    &              & \textbf{792} & \textbf{927.54} & \textbf{690}  \\ \bottomrule
			\end{tabular}
			\hspace{3mm}
			\begin{tabular}{@{}lrrrrr@{}}
				\toprule
				\textbf{}       &              & \textbf{$\Delta_{\text{OBJ}}$} & \textbf{$\Delta_{\text{CPU}}$} & \textbf{$\Delta_{\text{SOL}}$} \\ \cmidrule(r){1-2} \cmidrule(l){3-5}
				\textbf{SMALL}  & \textbf{ALL} & \textbf{-0.67} & \textbf{-15.14} & \textbf{0}   \\ \cmidrule(r){2-2} \cmidrule(l){3-5}
				& CM           & -0.63          & -39.60          & 0            \\
				& SM           & -0.70          & 146.97          & 0            \\
				& CI           & -0.71          & -24.36          & 0            \\
				& DI           & -0.62          & -3.09          & 0            \\
				& FI           & -1.00          & 45.79          & 0            \\
				& NI           & -0.28          & -39.40           & 0            \\
				\cmidrule(r){1-2} \cmidrule(l){3-5}
				\textbf{MEDIUM} & \textbf{ALL} & \textbf{0.94}  & \textbf{-85.54} & \textbf{9}   \\ \cmidrule(r){2-2} \cmidrule(l){3-5}
				& CM           & 1.07           & -86.65          & 0            \\
				& SM           & 0.79           & -84.03          & 9            \\
				& CI           & 1.05           & -86.33          & 5            \\
				& DI           & 0.83           & -84.63          & 4            \\
				& FI           & 0.21           & -84.91          & 5            \\
				& NI           & 1.87           & -85.94          & 4            \\
				\cmidrule(r){1-2} \cmidrule(l){3-5}
				\textbf{LARGE}  & \textbf{ALL} & \textbf{6.28}  & \textbf{-59.81}  & \textbf{90}  \\ \cmidrule(r){2-2} \cmidrule(l){3-5}
				& CM           & 6.29           & -65.41            & 30           \\ 
				& SM           & 6.02           & -51.10           & 60           \\
				& CI           & 5.73           & -58.53           & 45           \\
				& DI           & 6.83           & -61.10           & 45           \\
				& FI           & 6.16           & -42.32           & 49           \\
				& NI           & 7.61           & -78.36          & 41           \\
				\cmidrule(r){1-2} \cmidrule(l){3-5}
				\textbf{ALL}    &              & \textbf{1.23}  & \textbf{-69.59}  & \textbf{99} \\ \bottomrule
			\end{tabular}
		\end{adjustbox}
		\caption{The CSH results (left) and the comparison with the CPA for the \topPL{} (right).}
		\label{table:CSH_CPA_PL}
	\end{table}
	
	The results reported in Table~\ref{table:CSH_CPA_PL} shows that the CSH is able to compute high quality 
	solutions for all the small and medium instances, on average.
	Contrary to the results for the \topP{} instances, the CSH is able to further improve the gaps on larger
	instances.
	Regarding the running time, the CSH is always faster than the CPA on all instances.
	Finally, the CSH (as the VND) is capable to find 99 additional solutions with respect to the CPA,
	90 of which are solutions for the larger instances.

	Tables~\ref{table:VND_CPA_PL} and~\ref{table:CSH_CPA_PL} highlight an even more pronounced trend with respect
	to that for the \topP{} in terms of number of feasible solutions. Indeed, the proposed algorithms found 60 
	additional feasible solutions for the SM instances compared to 30 for the CM instances.
	In this case again, no relevant difference emerge in the number of feasible solutions found by the VND
	and CSH for the CI and DI instances. The same logic also follows for the FI and NI instances.
	Therefore, we can argue that the SM instances, that is those instances with mandatory nodes far from 
	each other, seem easier to solve than the other instances, confirming the analysis provided in \citet{Guastalla2024}.
	
	
	Finally, we provide a direct comparison between the VND and the CSH algorithms. 
	Table~\ref{table:CSH_VND} reports such a comparison on the \topP{} and \topPL{}
	instances respectively on the left and on the right.
	The results are summarised in the ALL rows, also by instance size and generation methods.
	The columns $\Delta_{\text{OBJ}}$ and $\Delta_{\text{CPU}}$ measure the average
	relative gaps (in percentage) of the CSH with respect to the VND for the objective function and 
	for the running time, respectively. 
	The column $\Delta_{\text{SOL}}$ measures the absolute gap between the two algorithms 
	in terms of number of solutions found.
	We would remark that a positive value in the columns $\Delta_{\text{OBJ}}$ and 
	$\Delta_{\text{SOL}}$ means that the CSH is better than the VND. On the contrary, 
	a negative value in the column $\Delta_{\text{CPU}}$ means that the CSH is better 
	than the VND.
	
	\begin{table}[!ht]
		\centering
		\begin{adjustbox}{width=0.65\textwidth,center}
			\begin{tabular}{@{}lrrrr@{}}
				\toprule
				\textbf{}       &              & \textbf{$\Delta_{\text{OBJ}}$} & \textbf{$\Delta_{\text{CPU}}$} & \textbf{$\Delta_{\text{SOL}}$} \\ \cmidrule(r){1-2} \cmidrule(l){3-5}
				\textbf{SMALL}  & \textbf{ALL} & \textbf{4.44}                  & \textbf{2482.05}               & \textbf{0}   \\ \cmidrule(r){2-2} \cmidrule(l){3-5}
				& CM           & 5.15                           & 3738.70                       & 0            \\
				& SM           & 3.65                           & 1763.47                        & 0            \\
				& CI           & 4.32                           & 3336.06                       & 0            \\
				& DI           & 4.56                           & 1900.64                        & 0            \\
				& -            & -                                & -                                & -            \\
				& -            & -                                & -                                & -            \\ \cmidrule(r){1-2} \cmidrule(l){3-5}
				\textbf{MEDIUM} & \textbf{ALL} & \textbf{3.37}                  & \textbf{382.48}               & \textbf{0}   \\ \cmidrule(r){2-2} \cmidrule(l){3-5}
				& CM           & 2.68                           & 415.41                        & 0            \\
				& SM           & 4.13                           & 348.70                        & 0            \\
				& CI           & 2.81                           & 327.90                        & 0            \\
				& DI           & 3.97                           & 449.72                        & 0            \\
				& -            & -                                & -                                & -            \\
				& -            & -                                & -                                & -            \\ \cmidrule(r){1-2} \cmidrule(l){3-5}
				\textbf{LARGE}  & \textbf{ALL} & \textbf{2.91}                  & \textbf{107.48}                & \textbf{0}   \\ \cmidrule(r){2-2} \cmidrule(l){3-5}
				& CM           & 2.54                           & 83.81                         & 0            \\
				& SM           & 5.54                           & 145.53                         & 0            \\
				& CI           & 3.81                           & 109.96                         & 0            \\
				& DI           & 2.13                           & 104.87                        & 0            \\
				& -            & -                                & -                                & -            \\
				& -            & -                                & -                                & -            \\ \cmidrule(r){1-2} \cmidrule(l){3-5}
				\textbf{ALL}    &              & \textbf{3.52}                  & \textbf{137.96}                & \textbf{0}               \\ \bottomrule
			\end{tabular}
			\hspace{3mm}
			\begin{tabular}{@{}lrrrr@{}}
				\toprule
				\textbf{}       &              & \textbf{$\Delta_{\text{OBJ}}$} & \textbf{$\Delta_{\text{CPU}}$} & \textbf{$\Delta_{\text{SOL}}$} \\ \cmidrule(r){1-2} \cmidrule(l){3-5}
				\textbf{SMALL}  & \textbf{ALL} & \textbf{3.47}                  & \textbf{1926.75}               & \textbf{0}   \\ \cmidrule(r){2-2} \cmidrule(l){3-5}
				& CM           & 3.98                           & 3043.59                        & 0            \\
				& SM           & 2.90                           & 1186.11                       & 0            \\
				& CI           & 3.54                           & 2271.07                        & 0            \\
				& DI           & 3.41                           & 1665.27                        & 0            \\
				& FI           & 4.36                           & 1509.01                        & 0            \\
				& NI           & 2.47                           & 2597.47                        & 0            \\ \cmidrule(r){1-2} \cmidrule(l){3-5}
				\textbf{MEDIUM} & \textbf{ALL} & \textbf{3.08}                  & \textbf{256.98}                & \textbf{0}   \\ \cmidrule(r){2-2} \cmidrule(l){3-5}
				& CM           & 2.59                           & 260.15                         & 0            \\
				& SM           & 3.66                           & 253.43                         & 0            \\
				& CI           & 3.16                           & 251.74                         & 0            \\
				& DI           & 2.99                           & 262.58                         & 0            \\
				& FI           & 2.38                           & 153.17                         & 0            \\
				& NI           & 3.96                           & 397.03                        & 0            \\ \cmidrule(r){1-2} \cmidrule(l){3-5}
				\textbf{LARGE}  & \textbf{ALL} & \textbf{2.15}                  & \textbf{44.59}                & \textbf{0}   \\ \cmidrule(r){2-2} \cmidrule(l){3-5}
				& CM           & 2.11                           & 26.55                         & 0            \\
				& SM           & 3.53                           & 71.51                         & 0            \\
				& CI           & 2.61                           & 42.41                         & 0            \\
				& DI           & 1.68                           & 47.01                         & 0            \\
				& FI           & 2.13                           & 35.91                         & 0            \\
				& NI           & 2.31                           & 76.44                         & 0            \\ \cmidrule(r){1-2} \cmidrule(l){3-5}
				\textbf{ALL}    &              & \textbf{3.05}                  & \textbf{69.31}                & \textbf{0}     \\ \bottomrule
			\end{tabular}
		\end{adjustbox}
		\caption{VND and CSH comparison for \topP{} (left) and for the \topPL{} (right).}
		\label{table:CSH_VND}
	\end{table}
	
	The results reported in Table~\ref{table:CSH_VND} show that the CSH is effective than the VND. On the contrary, the VND results much more efficient 
	than the CSH.
	We oberve that the relative gaps ($\Delta_{\text{OBJ}}$ and $\Delta_{\text{CPU}}$)
	are slightly better for the CSH on the \topP{} instances. This means that the CSH results 
	to be more effective but also less efficient for the \topP{}.
	Furthermore, the $\Delta_{\text{SOL}}$ is always zero since both algorithms are capable to 
	compute the same number of feasible solutions. This provides strong evidence about the 
	solvability of those instances for which the CPA is unable either to find a feasible 
	solution or to establish its infeasibility.

	\subsection{Computational results on TOP instances}
	\label{ssec:results-top}
	
	In this section, we provide a computational comparison among the best algorithms for the TOP and our VND and CSH algorithms.
	The idea is to provide a further validation of our algorithms.
	It is worth noting that considering the TOP as a baseline problem means that we are not considering the three features
	characterising the \topST{} on which our algorithm have been developed.
	
	To the best of our knowledge, the three best heuristic algorithms for the TOP are  
	\begin{inparaenum}[(i)]
		\item
		the PSOiA \citep{DANG2013332},
		\item
		the PMA \citep{KE2016155},
		\item 
		the HALNS \citep{HAMMAMI2020105034}.
	\end{inparaenum}
	In our computational experiments, we use the classical benchmark instances introduced by \citet{CHAO1996464}
	and those introduced by \citet{DANG2013332}. This benchmark set is composed by instances whose number of nodes 
	ranges from 101 to 400, generated in a similar way to those ones belonging to the set of \citeauthor{CHAO1996464}. 
	Regarding the three competitors, we consider the best value achieved across all the runs. All the comparisons
	are made by pairs of VND or CSH algorithms and one competitor.
	In order to provide a fair comparison, the reported running time of the competitors have been rescaled with 
	respect to our CPU frequency (2.1 GHz).
	Furthermore, all the three competitors have been tested multiple times using different seeds to prove their 
	robustness since they have some random ingredients. By consequence, we further rescaled their average 
	running times in accordance with the number of runs for each instance ($10$ runs for the PSOiA and the PMA, 
	and $20$ for HALNS).
	Finally, we are aware that the comparison between the CSH and its competitors may not be entirely fair. 
	However, the aim of this analysis is primarily qualitative (despite the fact that we compare the algorithms 
	running times) as our algorithms were not designed for TOP. In addition, the CSH has the intrinsic quality 
	of being parallel by design, unlike its competitors, who would require some effort to achieve some form of 
	parallelism.

	Table~\ref{table:CSH_VND_TOP} reports the results of the comparison between our algorithms and the three 
	competitors on the sets 4, 5, 6 and 7 of \citeauthor{CHAO1996464}'s instances. We have not considered the
	sets 1, 2 and 3 since they are considered the easiest ones.
	The columns $\Delta_{\text{MAX}}$ and $\Delta_{\text{AVG}}$ show the average relative gaps (in percentage) on the objective 
	function with respect to the best and the average solution values of the competitors ($\Delta_{\text{AVG}}$
	is available only for PSOiA). The column $\Delta_{\text{CPU}}$ reports the running time average relative differences (in percentage).
	A positive value of $\Delta_{\text{MAX}}$ and $\Delta_{\text{AVG}}$, and a negative value of $\Delta_{\text{CPU}}$ 
	indicates that our algorithm is better than the competitor.
	
	\begin{table}[!ht]
		\centering
		\begin{adjustbox}{width=.65\textwidth,center}
			\begin{tabular}{@{}lrrrrrrrrrr@{}}
				\toprule
				&                                  & \multicolumn{3}{c}{\textbf{PSOiA}}                             & \multicolumn{3}{c}{\textbf{PMA}}                               & \multicolumn{3}{c}{\textbf{HALNS}}                             \\
				& \textbf{Set} & \textbf{$\Delta_{\text{MAX}}$} & \textbf{$\Delta_{\text{AVG}}$} & \textbf{$\Delta_{\text{CPU}}$} & \textbf{$\Delta_{\text{MAX}}$} & \textbf{$\Delta_{\text{AVG}}$} & \textbf{$\Delta_{\text{CPU}}$} &\textbf{$\Delta_{\text{MAX}}$} & \textbf{$\Delta_{\text{AVG}}$} & \textbf{$\Delta_{\text{CPU}}$} \\ \cmidrule(r){1-2} \cmidrule(lr){3-5} \cmidrule(lr){6-8} \cmidrule(l){9-11}
				\multicolumn{1}{c}{\textbf{VND}} & \textbf{4}                       & -0.71                            & -0.68        & -57.72       & -0.71                            & -            & -31.30       & -0.71                            & -            & -1.93        \\
				\multicolumn{1}{c}{\textbf{}}    & \textbf{5}                       & -0.20                            & -0.05       & -53.54       & -0.20                            & -            & -18.40       & -0.20                            & -            & -32.35       \\
				\multicolumn{1}{c}{\textbf{}}    & \textbf{6}                       & 0.00                             & 0.00         & -56.64       & 0.00                             & -            & -54.43       & 0.00                             & -            & -27.78       \\
				\multicolumn{1}{c}{\textbf{}}    & \textbf{7}                       & -0.62                            & -0.82        & -40.81       & -0.62                            & -            & -14.14       & -0.62                            & -            & -36.10       \\ \cmidrule(r){1-2} \cmidrule(lr){3-5} \cmidrule(lr){6-8} \cmidrule(l){9-11}
				\multicolumn{1}{c}{\textbf{CSH}} & \textbf{4}                       & -0.02                            & 0.08         & -35.90       & -0.02                            & -            & 4.15         & -0.02                            & -            & 48.68        \\
				& \textbf{5}                       & 0.00                             & 0.01         & -46.50       & 0.00                             & -            & -6.04        & 0.00                             & -            & -22.10       \\
				& \textbf{6}                       & 0.00                             & 0.00         & -3.10        & 0.00                             & -            & 1.83         & 0.00                             & -            & 61.39        \\
				& \textbf{7}                       & -0.05                            & -0.04        & -4.37        & -0.05                            & -            & 38.70        & -0.05                            & -            & 3.23         \\ \bottomrule
			\end{tabular}
		\end{adjustbox}
		\caption{VND and CSH comparison with the three competitors on \citeauthor{CHAO1996464}'s instances.}
		\label{table:CSH_VND_TOP}
	\end{table}
	
	As reported in Table~\ref{table:CSH_VND_TOP}, the CSH results more effective in terms of solution quality 
	but less efficient in terms of average running time compared to the VND, confirming the behaviour observed for the
	\topST{}.
	The comparison with the three competitors shows that that the VND and the CSH are competitive with state-of-the-art 
	heuristic methods for the TOP in terms of both solution quality and running time.
	Even if our algorithms do not improve on any best known solution from the literature (because many of them are probably 
	optimal objective function values), the results are, on average, very close to the ones reported by the competitors, 
	especially for the CSH.
	Finally, regarding the average running times, the VND seems to be the fastest solution method across all sets and 
	algorithms, while the CSH seems to be more efficient than the PSOiA but slightly slower than the PMA and HALNS on 
	the average.

	Table~\ref{table:VND_TOP_LARGE} reports the results of the comparison between our algorithms and the three 
	competitors on the \citeauthor{DANG2013332}'s instances clustering them by size (number of nodes).
	The table has the same structure as Table~\ref{table:CSH_VND_TOP}.
	
	\begin{table}[!ht]
		\begin{adjustbox}{width=0.7\textwidth,center}
			\begin{tabular}{@{}llrrrrrrrrr@{}}
				\toprule
				\multicolumn{1}{l}{}          &              & \multicolumn{3}{c}{\textbf{PSOiA}}                                                       & \multicolumn{3}{c}{\textbf{PMA}}                                                         & \multicolumn{3}{c}{\textbf{HALNS}}                                                      \\
				& \textbf{Size} & {\textbf{$\Delta_{\text{MAX}}$}} & \textbf{$\Delta_{\text{AVG}}$} & \textbf{$\Delta_{\text{CPU}}$}          & {\textbf{$\Delta_{\text{MAX}}$}} & \textbf{$\Delta_{\text{AVG}}$} & \textbf{$\Delta_{\text{CPU}}$}              & {\textbf{$\Delta_{\text{MAX}}$}} & \textbf{$\Delta_{\text{AVG}}$} & \textbf{$\Delta_{\text{CPU}}$}            \\ \cmidrule(r){1-2} \cmidrule(lr){3-5} \cmidrule(lr){6-8} \cmidrule(l){9-11}
				\textbf{VND}                  & 101 to 200   & -3.22 & -2.91 & -34.43 & -3.23 & -2.99 & 27.87  & -3.23 & -3.10 & 19.24 \\
				\textbf{}                     & 201 to 300   & -3.25 & -2.90 & 19.40  & -3.26 & -2.89 & 213.50 & -3.26 & -3.05 & 53.24 \\
				\textbf{}                     & 301 to 400   & -6.94 & -6.02 & -82.49 & -7.11 & -6.26 & 45.47 & -7.11 & -6.34 & -8.76 \\ \cmidrule(r){1-2} \cmidrule(lr){3-5} \cmidrule(lr){6-8} \cmidrule(l){9-11}
				\textbf{CSH}                  & 101 to 200   & -0.41                     & -0.08                     & 49.56                      & -0.42                     & -0.17                    & 167.37                     & -0.42                     & -0.28                     & 138.77                    \\
				\textbf{}                     & 201 to 300   & -1.00                    & -0.64                     & 135.83                     & -1.02                     & -0.63                     & 343.35                     & -1.02                     & -0.80                     & 145.21                    \\
				\multicolumn{1}{l}{\textbf{}} & 301 to 400   & -2.52                     & -1.54                    & -74.54                     & -2.70                     & -1.80                     & 86.75                      & -2.70                     & -1.88                     & 28.44                     \\ \bottomrule
			\end{tabular}
		\end{adjustbox}
		\caption{VND and CSH comparison with the three competitors on \citeauthor{DANG2013332}'s instances.}
		\label{table:VND_TOP_LARGE}
	\end{table}
	
	As shown in Table~\ref{table:VND_TOP_LARGE}, the VND is generally less effective and efficient than the 
	other approaches, especially in the last subset (it is faster only compared to the PSOiA in the first 
	and last sets, on average). 
	In contrast, the CSH remains competitive with the state-of-the-art algorithms, although it is, on average, 
	slower than the other approaches (it is slightly faster than PSOiA only in the last set). Actually, the 
	worst decrease in the objective function among the different subsets is around 2.70\% taking into account 
	the $\Delta_{\text{MAX}}$ or 1.88\% taking into account the $\Delta_{\text{AVG}}$.

	
	\section{Conclusions}
	\label{sec:end}
	
	In this work, we introduced two different algorithms for the \topST{} which exploit a heuristic based on the 
	mandatory nodes allocation to compute an initial solution. This algorithm is able to quickly find a feasible 
	solution even if the feasibility problem of the TOP-ST-MIN has been proved to be NP-complete. 
	The first algorithm is a metaheuristic approach based on the Variable Neighbourhood Descent (VND) methodology.
	The VND exploits two families of neighbourhoods (profit-based and cost-based) to explore feasible and infeasible
	regions that use a couple of variable-size tabu lists to fix some nodes inside or outside the solution.
	The second one is a matheuristic approach that we called Cuts-Separation Heuristic (CSH) that generates a
	sequence of restricted problems in which adding subtours elimination constraints guides the algorithm 
	towards convergence to an optimal solution. Each solution in the sequence is then improved by the VND.
	We developed also a multi-thread version of the algorithm which exploits its intrinsic quality of being 
	parallel by design.
	Finally, we give an extensive comparison between the two algorithms and the Cutting-Plane Algorithm (CPA) 
	reported in \citet{Guastalla2024}. Furthermore, we provide an extensive comparison of the proposed algorithms 
	adapted to solve the TOP with the three best algorithms in literature for the problem.

	The results clearly show the complementarity of both approaches: the VND is much efficient but less effective 
	with respect to the CSH, i.e., the CSH usually finds better quality solutions but requires larger computational 
	times. Indeed, the VND is always faster than the CPA for both problems \topP{} and \topPL{}, and for every instance 
	size. 
	On the contrary, the CSH is as effective as the CPA for the \topP{} but results to be much better for the \topPL{}, 
	especially for the large instances.
	Finally, our heuristic approaches are able to find 48 additional feasible solutions for the \topP{} and 99 for 
	the \topPL{} with respect to the CPA proving their effectiveness at finding feasible solutions, even for larger 
	and more difficult instances. 
	The CSH is also competitive with the state-of-art TOP heuristic algorithms, in the sense that it provides solutions 
	with objective values within a few percent of the best known results in the literature. On the contrary, the VND 
	turns out to be competitive only on the instances from \citet{CHAO1996464}.

	As discussed in \citet{Nicosia2017}, modelling fairness is in itself challenging since the concept of fairness
	can vary in accordance with the context of the problem. This results in modelling approaches providing solutions 
	with different quality. Work is under way to exploit the flexibility of our algorithms to
	investigate the impact on the solutions of a fairness perspective in order to evaluate the so called price of
	fairness \citep{Bertsimas2011}.

	%
	
	\bibliographystyle{itor}
	\bibliography{references}

\begin{thebibliography}{47}
\expandafter\ifx\csname natexlab\endcsname\relax\def\natexlab#1{#1}\fi
\expandafter\ifx\csname url\endcsname\relax
  \def\url#1{\texttt{#1}}\fi
\expandafter\ifx\csname urlprefix\endcsname\relax\def\urlprefix{URL }\fi
\providecommand{\eprint}[2][]{\url{#2}}
\bibitem[{Aghezzaf and Fahim(2015)}]{Aghezzaf2015}
Aghezzaf, B., Fahim, H., 2015.
\newblock Solving the capacitated team orienteering problem with time windows
  through variable neighborhood search.
\newblock \textit{International Review on Computers and Software (IRECOS)} 10,
  1134.
\bibitem[{Aringhieri et~al.(2024)Aringhieri, Bigharaz, Druetto, Duma, Grosso
  and Guastalla}]{DSTCP}
Aringhieri, R., Bigharaz, S., Druetto, A., Duma, D., Grosso, A., Guastalla, A.,
  2024.
\newblock The daily swab test collection problem.
\newblock \textit{Annals of Operations Research} 335, 1449--1470.
\bibitem[{Aringhieri et~al.(2022a)Aringhieri, Bigharaz, Duma and
  Guastalla}]{ARP}
Aringhieri, R., Bigharaz, S., Duma, D., Guastalla, A., 2022a.
\newblock Fairness in ambulance routing for post disaster management.
\newblock \textit{Central European Journal of Operations Research} 30,
  189--211.
\bibitem[{Aringhieri et~al.(2022b)Aringhieri, Bigharaz, Duma and
  Guastalla}]{ODS}
Aringhieri, R., Bigharaz, S., Duma, D., Guastalla, A., 2022b.
\newblock Novel applications of the team orienteering problem in health care
  logistics.
\newblock In Amorosi, L., Dell'Olmo, P. and Lari, I. (eds),
  \textit{Optimization in Artificial Intelligence and Data Sciences}, Springer
  International Publishing, Cham, pp. 235--245.
\bibitem[{Aringhieri et~al.(2015)Aringhieri, Landa, Soriano, T\`anfani and
  Testi}]{ORP-COR}
Aringhieri, R., Landa, P., Soriano, P., T\`anfani, E., Testi, A., 2015.
\newblock A two level {M}etaheuristic for the {O}perating {R}oom {S}cheduling
  and {A}ssignment {P}roblem.
\newblock \textit{Computers \& Operations Research} 54, 21--34.
\bibitem[{Bertsimas et~al.(2011)Bertsimas, Farias and
  Trichakis}]{Bertsimas2011}
Bertsimas, D., Farias, V.F., Trichakis, N., 2011.
\newblock The price of fairness.
\newblock \textit{Operations Research} 59, 1, 17--31.
\bibitem[{Chao et~al.(1996)Chao, Golden and Wasil}]{CHAO1996464}
Chao, I.M., Golden, B.L., Wasil, E.A., 1996.
\newblock The team orienteering problem.
\newblock \textit{European Journal of Operational Research} 88, 464--474.
\bibitem[{Corrêa et~al.(2024)Corrêa, Dong, Iori, Santos, Yagiura and
  Zucchi}]{Correa2024}
Corrêa, V.H., Dong, H., Iori, M., Santos, A., Yagiura, M., Zucchi, G., 2024.
\newblock An iterated local search for a multi‐period orienteering problem
  arising in a car patrolling application.
\newblock \textit{Networks} 83, 153--168.
\bibitem[{Dang et~al.(2011)Dang, Guibadj and Moukrim}]{Dang2011}
Dang, D.C., Guibadj, R.N., Moukrim, A., 2011.
\newblock A pso-based memetic algorithm for the team orienteering problem.
\newblock In \textit{Applications of Evolutionary Computation}, Springer Berlin
  Heidelberg, pp. 471--480.
\bibitem[{Dang et~al.(2013)Dang, Guibadj and Moukrim}]{DANG2013332}
Dang, D.C., Guibadj, R.N., Moukrim, A., 2013.
\newblock An effective pso-inspired algorithm for the team orienteering
  problem.
\newblock \textit{European Journal of Operational Research} 229, 2, 332--344.
\bibitem[{Dantzig et~al.(1954)Dantzig, Fulkerson and Johnson}]{dfj1954}
Dantzig, G., Fulkerson, R., Johnson, S., 1954.
\newblock Solution of a large-scale traveling-salesman problem.
\newblock \textit{Journal of the Operations Research Society of America} 2,
  393--410.
\bibitem[{Duarte et~al.(2018)Duarte, S{\'a}nchez-Oro, Mladenovi{\'{c}} and
  Todosijevi{\'{c}}}]{Duarte2018}
Duarte, A., S{\'a}nchez-Oro, J., Mladenovi{\'{c}}, N., Todosijevi{\'{c}}, R.,
  2018.
\newblock Variable Neighborhood Descent, Springer Cham.
\newblock pp. 341--367.
\bibitem[{Dutta et~al.(2020)Dutta, Barma, Mukherjee, Kar and De}]{Dutta2020}
Dutta, J., Barma, P., Mukherjee, A., Kar, S., De, T., 2020.
\newblock A multi-objective open set orienteering problem.
\newblock \textit{Neural Computing and Applications} 32, 13953–--13969.
\bibitem[{Expósito~Márquez et~al.(2019)Expósito~Márquez, Mancini, Brito and
  Moreno-Pérez}]{Exposito2019}
Expósito~Márquez, A., Mancini, S., Brito, J., Moreno-Pérez, J., 2019.
\newblock A fuzzy grasp for the tourist trip design with clustered pois.
\newblock \textit{Expert Systems with Applications} 127.
\bibitem[{Feo~Flushing et~al.(2021)Feo~Flushing, Gambardella and
  Di~Caro}]{Feo2021}
Feo~Flushing, E., Gambardella, L., Di~Caro, G., 2021.
\newblock Spatially-distributed missions with heterogeneous multi-robot teams.
\newblock \textit{IEEE Access} 9, 67327--67348.
\bibitem[{Guastalla et~al.(2025)Guastalla, Aringhieri and
  Hosteins}]{Guastalla2024}
Guastalla, A., Aringhieri, R., Hosteins, P., 2025.
\newblock The team orienteering problem with service times and mandatory \&
  incompatible nodes.
\newblock \textit{Computers \& Operations Research} 183, 107170.
\bibitem[{G\"{u}ndling and Witzel(2020)}]{Gndling2020TimeDependentTT}
G\"{u}ndling, F., Witzel, T., 2020.
\newblock {Time-Dependent Tourist Tour Planning with Adjustable Profits}.
\newblock In \textit{20th Symposium on Algorithmic Approaches for
  Transportation Modelling, Optimization, and Systems},  \textit{Open Access
  Series in Informatics}. Vol.~85, pp. 1--14.
\bibitem[{Hammami et~al.(2020)Hammami, Rekik and Coelho}]{HAMMAMI2020105034}
Hammami, F., Rekik, M., Coelho, L.C., 2020.
\newblock A hybrid adaptive large neighborhood search heuristic for the team
  orienteering problem.
\newblock \textit{Computers \& Operations Research} 123, 105034.
\bibitem[{Hanafi et~al.(2020)Hanafi, Mansini and Zanotti}]{Hanafi2020}
Hanafi, S., Mansini, R., Zanotti, R., 2020.
\newblock The multi-visit team orienteering problem with precedence
  constraints.
\newblock \textit{European Journal of Operational Research} 282, 515--529.
\bibitem[{Heris et~al.(2021)Heris, Ghannadpour, Bagheri and
  Zandieh}]{Heris2021}
Heris, F., Ghannadpour, S., Bagheri, M., Zandieh, F., 2021.
\newblock A new accessibility based team orienteering approach for urban
  tourism routes optimization (a real life case).
\newblock \textit{Computers \& Operations Research} 138, 105620.
\bibitem[{Jin and Thomas(2019)}]{Jin2019}
Jin, H., Thomas, B., 2019.
\newblock Team orienteering with uncertain rewards and service times with an
  application to phlebotomist intrahospital routing.
\newblock \textit{Networks} 73, 453--465.
\bibitem[{Ke et~al.(2016)Ke, Zhai, Li and Chan}]{KE2016155}
Ke, L., Zhai, L., Li, J., Chan, F.T., 2016.
\newblock Pareto mimic algorithm: An approach to the team orienteering problem.
\newblock \textit{Omega} 61, 155--166.
\bibitem[{Kim et~al.(2013)Kim, Li and Johnson}]{KIM20133065}
Kim, B.I., Li, H., Johnson, A.L., 2013.
\newblock An augmented large neighborhood search method for solving the team
  orienteering problem.
\newblock \textit{Expert Systems with Applications} 40, 8, 3065--3072.
\bibitem[{Kotiloglu et~al.(2017)Kotiloglu, Lappas, Pelechrinis and
  Repoussis}]{Kotiloglu2017}
Kotiloglu, S., Lappas, T., Pelechrinis, K., Repoussis, P., 2017.
\newblock Personalized multi-period tour recommendations.
\newblock \textit{Tourism Management} 62, 76--88.
\bibitem[{Li et~al.(2024)Li, Zhu, Peng, Wang, Zhen and
  Demeulemeester}]{LI2024793}
Li, J., Zhu, J., Peng, G., Wang, J., Zhen, L., Demeulemeester, E., 2024.
\newblock Branch-price-and-cut algorithms for the team orienteering problem
  with interval-varying profits.
\newblock \textit{European Journal of Operational Research} 319, 3, 793--807.
\bibitem[{Liao and Zheng(2018)}]{Liao2018}
Liao, Z., Zheng, W., 2018.
\newblock Using a heuristic algorithm to design a personalized day tour route
  in a time-dependent stochastic environment.
\newblock \textit{Tourism Management} 68, 284--300.
\bibitem[{Lim et~al.(2018)Lim, Chan, Leckie and Karunasekera}]{Lim2018}
Lim, K.H., Chan, J., Leckie, C., Karunasekera, S., 2018.
\newblock Personalized trip recommendation for tourists based on user
  interests, points of interest visit durations and visit recency.
\newblock \textit{Knowledge and Information Systems} 54, 375--406.
\bibitem[{Lin and Yu(2017)}]{LIN2017195}
Lin, S.W., Yu, V.F., 2017.
\newblock Solving the team orienteering problem with time windows and mandatory
  visits by multi-start simulated annealing.
\newblock \textit{Computers \& Industrial Engineering} 114, 195--205.
\bibitem[{Lu et~al.(2018)Lu, Benlic and Wu}]{LU201854}
Lu, Y., Benlic, U., Wu, Q., 2018.
\newblock A memetic algorithm for the orienteering problem with mandatory
  visits and exclusionary constraints.
\newblock \textit{European Journal of Operational Research} 268, 1, 54--69.
\bibitem[{Nicosia et~al.(2017)Nicosia, Pacifici and Pferschy}]{Nicosia2017}
Nicosia, G., Pacifici, A., Pferschy, U., 2017.
\newblock Price of fairness for allocating a bounded resource.
\newblock \textit{European Journal of Operational Research} 257, 3, 933--943.
\bibitem[{Palomo-Martínez et~al.(2015)Palomo-Martínez, Salazar-Aguilar,
  Laporte and Langevin}]{Palomo2015}
Palomo-Martínez, P., Salazar-Aguilar, A., Laporte, G., Langevin, A., 2015.
\newblock A hybrid variable neighborhood search for the orienteering problem
  with mandatory visits and exclusionary constraints.
\newblock \textit{Computers \& Operations Research} 78, 408--419.
\bibitem[{Reyes-Rubiano et~al.(2018)Reyes-Rubiano, Ospina-Trujillo, Faulin,
  Mozos, Panadero and Juan}]{Reyes2018}
Reyes-Rubiano, L., Ospina-Trujillo, C., Faulin, J., Mozos, J., Panadero, J.,
  Juan, A., 2018.
\newblock The team orienteering problem with stochastic service times and
  driving-range limitations: a simheuristic approach.
\newblock In \textit{2018 Winter Simulation Conference (WSC)}, IEEE, pp.
  3025--3035.
\bibitem[{Ruiz-Meza et~al.(2021)Ruiz-Meza, Brito and Montoya-Torres}]{Ruiz2021}
Ruiz-Meza, J., Brito, J., Montoya-Torres, J., 2021.
\newblock A grasp to solve the multi-constraints multi-modal team orienteering
  problem with time windows for groups with heterogeneous preferences.
\newblock \textit{Computers \& Industrial Engineering} 162, 107776.
\bibitem[{Shen(1965)}]{Lin1965}
Shen, L., 1965.
\newblock Computer solutions of the traveling salesman problem.
\newblock \textit{The Bell System Technical Journal} 44, 10, 2245--2269.
\bibitem[{Stavropoulou et~al.(2019)Stavropoulou, Repoussis and
  Tarantilis}]{Stavropoulou2019}
Stavropoulou, F., Repoussis, P., Tarantilis, C., 2019.
\newblock The vehicle routing problem with profits and consistency constraints.
\newblock \textit{European Journal of Operational Research} 274, 340--356.
\bibitem[{Sun et~al.(2018)Sun, Veelenturf, Dabia and {Van
  Woensel}}]{SUN20181058}
Sun, P., Veelenturf, L.P., Dabia, S., {Van Woensel}, T., 2018.
\newblock The time-dependent capacitated profitable tour problem with time
  windows and precedence constraints.
\newblock \textit{European Journal of Operational Research} 264, 3, 1058--1073.
\bibitem[{Sylejmani et~al.(2017)Sylejmani, Dorn and Musliu}]{Sylejmani2017}
Sylejmani, K., Dorn, J., Musliu, N., 2017.
\newblock Planning the trip itinerary for tourist groups.
\newblock \textit{Information Technology \& Tourism} 17.
\bibitem[{Thompson and Galeazzi(2019)}]{Thompson2019}
Thompson, F., Galeazzi, R., 2019.
\newblock Robust mission planning for autonomous marine vehicle fleets.
\newblock \textit{Robotics and Autonomous Systems} 124, 103404.
\bibitem[{Vansteenwegen and Gunawan(2019)}]{OP2019}
Vansteenwegen, P., Gunawan, A., 2019.
\newblock \textit{Orienteering Problems: Models and Algorithms for Vehicle
  Routing Problems with Profits}.
\newblock Springer Nature Switzerland AG.
\bibitem[{Vansteenwegen and Souffriau(2011)}]{Vansteenwegen2011}
Vansteenwegen, P., Souffriau, W., 2011.
\newblock The orienteering problem: A survey.
\newblock \textit{European Journal of Operational Research} 209, 1--10.
\bibitem[{Vidal et~al.(2014)Vidal, Crainic, Gendreau and Prins}]{VIDAL2014658}
Vidal, T., Crainic, T.G., Gendreau, M., Prins, C., 2014.
\newblock A unified solution framework for multi-attribute vehicle routing
  problems.
\newblock \textit{European Journal of Operational Research} 234, 3, 658--673.
\bibitem[{Wang et~al.(2017)Wang, Lau and Cheng}]{Wang2017}
Wang, W., Lau, H.C., Cheng, S.F., 2017.
\newblock Exact and heuristic approaches for the multi-agent orienteering
  problem with capacity constraints.
\newblock In \textit{2017 IEEE Symposium Series on Computational Intelligence
  (SSCI)}, pp. 1--7.
\bibitem[{Wisittipanich and Boonya(2020)}]{Wisittipanich2020}
Wisittipanich, W., Boonya, C., 2020.
\newblock Multi-objective tourist trip design problem in chiang mai city.
\newblock \textit{IOP Conference Series: Materials Science and Engineering}
  895, 012014.
\bibitem[{Yu et~al.(2019)Yu, Fang, Zhu and Ma}]{YU2019488}
Yu, Q., Fang, K., Zhu, N., Ma, S., 2019.
\newblock A matheuristic approach to the orienteering problem with service time
  dependent profits.
\newblock \textit{European Journal of Operational Research} 273, 2, 488--503.
\bibitem[{Yu et~al.(2024)Yu, Salsabila, Lin and Gunawan}]{YU2024121996}
Yu, V.F., Salsabila, N.Y., Lin, S.W., Gunawan, A., 2024.
\newblock Simulated annealing with reinforcement learning for the set team
  orienteering problem with time windows.
\newblock \textit{Expert Systems with Applications} 238, 121996.
\bibitem[{Zheng et~al.(2020)Zheng, Ji, Lin, Wang and Yu}]{Zheng2020}
Zheng, W., Ji, H., Lin, C., Wang, W., Yu, B., 2020.
\newblock Using a heuristic approach to design personalized urban tourism
  itineraries with hotel selection.
\newblock \textit{Tourism Management} 76, 10956.
\bibitem[{Zheng et~al.(2017)Zheng, Liao and Qin}]{ZHENG2017335}
Zheng, W., Liao, Z., Qin, J., 2017.
\newblock Using a four-step heuristic algorithm to design personalized day tour
  route within a tourist attraction.
\newblock \textit{Tourism Management} 62, 335--349.

\end{thebibliography}

	
	
\end{document}